\documentclass{article}

\usepackage{comment}
\usepackage{fullpage}
\usepackage{aageneral}
\usepackage{aamath}
\usepackage{mhequ}
\usepackage{aachem}
\usepackage{authblk}

\def\dpr{{\delta_0'}}
\def\PPG{\PP_t^G}
\newcommand{\cf}[1]{{\iffalse #1 \fi}}
\newcommand{\obar}{\overline}

\mathtoolsset{showonlyrefs=true}

\title{\Large Seemingly stable chemical kinetics can be stable, marginally stable, or unstable
}
\author{Andrea~Agazzi}
\author{Jonathan C. Mattingly}

\affil{Mathematics Department, Duke University}

\date{\today}

\begin{document}
\maketitle

\begin{abstract}
  We present three examples of chemical reaction networks whose ordinary differential equation scaling limit are almost identical and in all cases stable. Nevertheless, the Markov jump processes associated to these reaction networks display the full range of behaviors: one is stable (positive recurrent), one is unstable (transient) and one is marginally stable (null recurrent).
 We study these differences and characterize the invariant measures by Lyapunov function techniques. In particular, we design a natural set of such functions which scale homogeneously to infinity, taking advantage of the same scaling behavior of the reaction rates.
\end{abstract}


\section{Introduction}\label{s:networks}

The recent progress in experimental biology and bioinformatics has sparked renewed interest the theoretical description of cellular biochemical processes. A common approach in this sense is to model the dynamics of chemical reactions through mass action kinetics. The most popular formulation of this model is deterministic and describes the dynamics of the network through a set of ordinary differential equations. These systems of equations approximate the interactions of the individual molecules involved in the reaction network \cite{either86}. An alternative mass action kinetics framework takes into account the discrete nature of chemical systems by representing their state as the number of molecules of each species that are present in the reactor. In this formulation, transitions occur when molecules combine as inputs of a reaction and transform into its outputs, resulting in a jump in the state of the system. The dynamics of these discrete systems can be modeled stochastically as jump Markov processes \cite[Sec.~11, Example~C]{either86} whose jumping rates are specified, under the mass action kinetics assumption, by the structure of the chemical network being modeled. It is well know that, for large number of molecules, the dynamics of stochastic mass action kinetics converge to their deterministic counterpart \cite{agazzi17, either86}. However, when the number of molecules is finite the dynamics of the two families of models can be qualitatively different. \aa{However, networks with identical asymptotic behavior in the continuum approximation can display radically different dynamics when the number of molecules is finite.} This is the object of study of this paper.

The study of dynamical properties of mass action systems was greatly advanced with the work of Horn and Jackson \cite{horn72} and Feinberg \cite{feinberg87}. Asymptotic studies in the stochastic setting appear in the probability literature with the work of Kurtz \cite{either86}. More recently, results have appeared on the existence and characterization of the invariant measure \cite{anderson10, cappelletti16} of some classes of  Chemical Reaction Networks, abbreviated henceforth as \abbr{crn}s. Furthermore, in \cite{kim17}, criteria for recurrence of stochastic mass action models based on the geometry of the underlying \abbr{crn} have been developed, working toward the proof of the \emph{recurrence conjecture} \cite{kim17}. In this paper we study a specific family of examples displaying radically different recurrence properties in the discrete framework despite their stable and almost identical behaviour when modeled deterministically. The different asymptotic behavior stems from the different behavior in the neighborhood of the horizontal axis which the deterministic dynamics, obtained through a scaling procedure, ignores. This behavior is similar to the existence of a boundary-layer in singular perturbation theory. We note that a similar example has independently been presented in the recently submitted paper \cite{anderson18}. We then study the invariant measure of the presented networks and their convergence to it. In doing so we develop a method for the construction of Foster-Lyapunov functions (henceforth referred to as simply Lyapunov functions) that to the best of the knowledge of the authors has never been applied for the study of the stability of this class of systems. This method is based on dynamic principles and capable of producing Lyapunov functions which address delicate boundary cases between stability and instability.

\subsection{The model}

We consider a set of $d$ \emph{species} $\S := \{s_1,\dots, s_d\}$ whose interactions are described by a set of $m$ \emph{reactions} $\R := \{r_1, \dots, r_m\}$. Each chemical reaction $r \in \R$ describes how molecules of different species in the reactor combine as inputs of $r$ to form its outputs. Then, letting $\Nn_0$ denote the set of nonnegative integers, we represent any $r \in \R$ as
\begin{equ}
  r := \pg{\sum_{i = 1}^d (c_\inn^r)_i \rightharpoonup \sum_{i = 1}^d (c_\outt^r)_i}\,,
\end{equ}
where $c_\inn^r, c_\outt^r \in \Nn_0^d$ count the multiplicity of each species as input and output of the reaction $r$ and are respectively referred to as the \emph{input} and \emph{output complex} of $r$. We also denote by $\C := \{c_\#^r~:~r \in \R,\,\# \in \{\inn, \outt\}\}$ the \emph{set of complexes}. Finally, we uniquely identify a \abbr{crn} with the corresponding triplet \net.

The main object of study of this paper is the behavior of the process $X_t \in \Nn_0^d$, counting in each of its components the number of molecules in the corresponding species at time $t$. The effect of a reaction $r \in \R$ on the state of the network is encoded by the respective \emph{reaction vector} $c^r := r_\outt^r - c_\inn^r$\,: When that reaction occurs, the state of the network jumps as $X_t \to X_t + c^r$. The probabilistic dynamics of $X_t$ is modeled as a  jump Markov process. The generator of this process is given under the mass action kinetics assumption by
\begin{equs}\label{e:generator}
 \LL f(x) &= \sum_{r \in \R} \Lambda_r(x) \pc{f\pc{x + c_r} - f(x)} = \sum_{r \in \R} \Lambda_r(x) \Delta_r f(x)\,,
\end{equs}
for a function $f\colon\mathbb N_0^d \to \Rr$, and \emph{reaction rates} $\{\Lambda_r(N_t)\}_r$ defined by
\begin{equ}
 \Lambda_r(N_t) = \kappa_r \prod_{i = 1}^d \binom{(N_t)_i }{ (c^r_{\rm in})_i} (c^r_{\rm in})_i! ~,
 \label{e:rates}
\end{equ}
for \emph{reaction rate constants} $\kappa_r \in \Rr_{>0}$ and where $\binom{a}{b}$ is the binomial coefficient, set to $0$ if $b \not \in [0,a]$. Accordingly, we define the Markov transition kernel $\PP_t$ associated to the process $X_t$. The left-action of $\PP_t$ as a linear operator on the space of signed measures on $\RR^d$ and  right-action of $\PP_t$ on the space of bounded measurable functions will be denoted by
\begin{equ}
  (\mu\PP_t ) (A) = \int_{\Nn_0^d} \PP_t(x, A) \mu(\dd x) \qquad \text{and} \qquad   (\PP_t f) (x) = \int_{\Nn_0^d} f(y) \PP_t(x, \dd y)\,,
\end{equ}
for any measure $\mu$, set $A$ and bounded measurable function $f$.

The generator \eref{e:generator} and the rates \eref{e:rates} can be rescaled in a natural way \cite{either86} to describe the dynamics of the concentration vector $v^{{-1}} X_t \in (v^{-1} \Nn_0)^d$ in a reactor of volume $v > 0$. Denoting throughout by $\RR$ the set of nonnegative real numbers, it is well known that in the limit $v \to \infty$  the sample paths of $X_t^v$ with initial condition $\lim_{v \to \infty} X_0^v = x_0 \in \RR^d$ converge through a functional law of large numbers to the deterministic trajectories $x(t) \in \RR^d$ of the system of ordinary differential equations
\begin{equ}
\frac{\dd x}{\dt}  = \sum_{r \in \mathcal \R} \lambda_r(x) {c}^r \,, \qquad \text{with }  \quad\lambda_r(x) := \kappa_r \prod_{i = 1}^d {x}_i^{(c^r_\inn)_i} \quad\text{where }x=(x_1,\dots,x_d),
\label{e:ma}
\end{equ}
and initial condition $x(0) = x_0$, provided that a solution to \eref{e:ma} exists for the time interval of interest \cite[§11, Thm. 2.1]{either86}. \eref{e:ma} are referred to as \emph{deterministic mass action kinetics equations} and the regime $v \to \infty$ as the \emph{fluid limit} \cite{either86}.

\subsection{The examples and main results}

%

\begin{figure}[h!]
 \centering
  \def\svgwidth{.5\textwidth}
  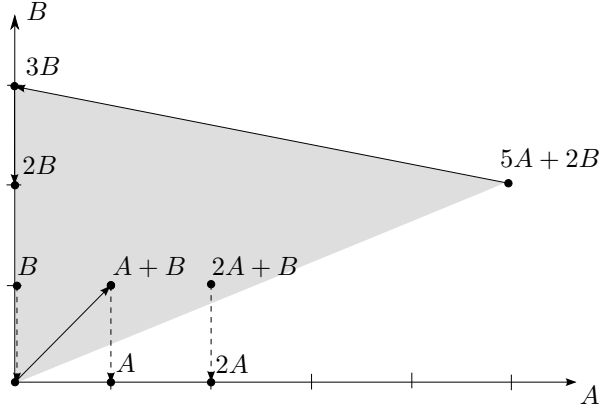
 \caption{The complex graph of the networks in \eref{e:e0}--\eref{e:e2}. The common reactions in the networks are displayed as solid arrows, while the reactions that are different in the three examples are dashed. The existence of the reaction $\emptyset \rightharpoonup A+B$ makes the networks asiphonic, and the fact that all the arrows starting on the faces of the reaction polytope (the grey triangle) point inwards makes the network Strongly Endotactic. }
 \label{f:f1}
\end{figure}

In this paper we consider the following \abbr{crn}s:
\begin{equs}
 \text{\abbr{CRN0}}: \qquad \emptyset \rightharpoonup A + B\,, \qquad &B \rightharpoonup \emptyset \,,  \qquad & 5A + 2B \rightharpoonup 3B \rightharpoonup 2A~, \label{e:e0}\\
 \text{\abbr{CRN1}}: \qquad \emptyset \rightharpoonup A + B\,, \qquad &A+B \rightharpoonup A\,, \qquad & 5A + 2B \rightharpoonup 3B \rightharpoonup 2A~, \label{e:e1}\\
 \text{\abbr{CRN2}}: \qquad \emptyset \rightharpoonup A + B\,, \qquad &2A+B \rightharpoonup 2A\,,  \qquad & 5A + 2B \rightharpoonup 3B \rightharpoonup 2A~. \label{e:e2}\\
\end{equs}
For each network, we number the reactions from left to right, obtaining $\R = \{r_1,r_2,r_3,r_4\}$. The networks are displayed in \fref{f:f1}.

\begin{figure}
 \centering
 \begin{subfigure}{.32\textwidth}
  \centering
  \includegraphics[width=.8\linewidth]{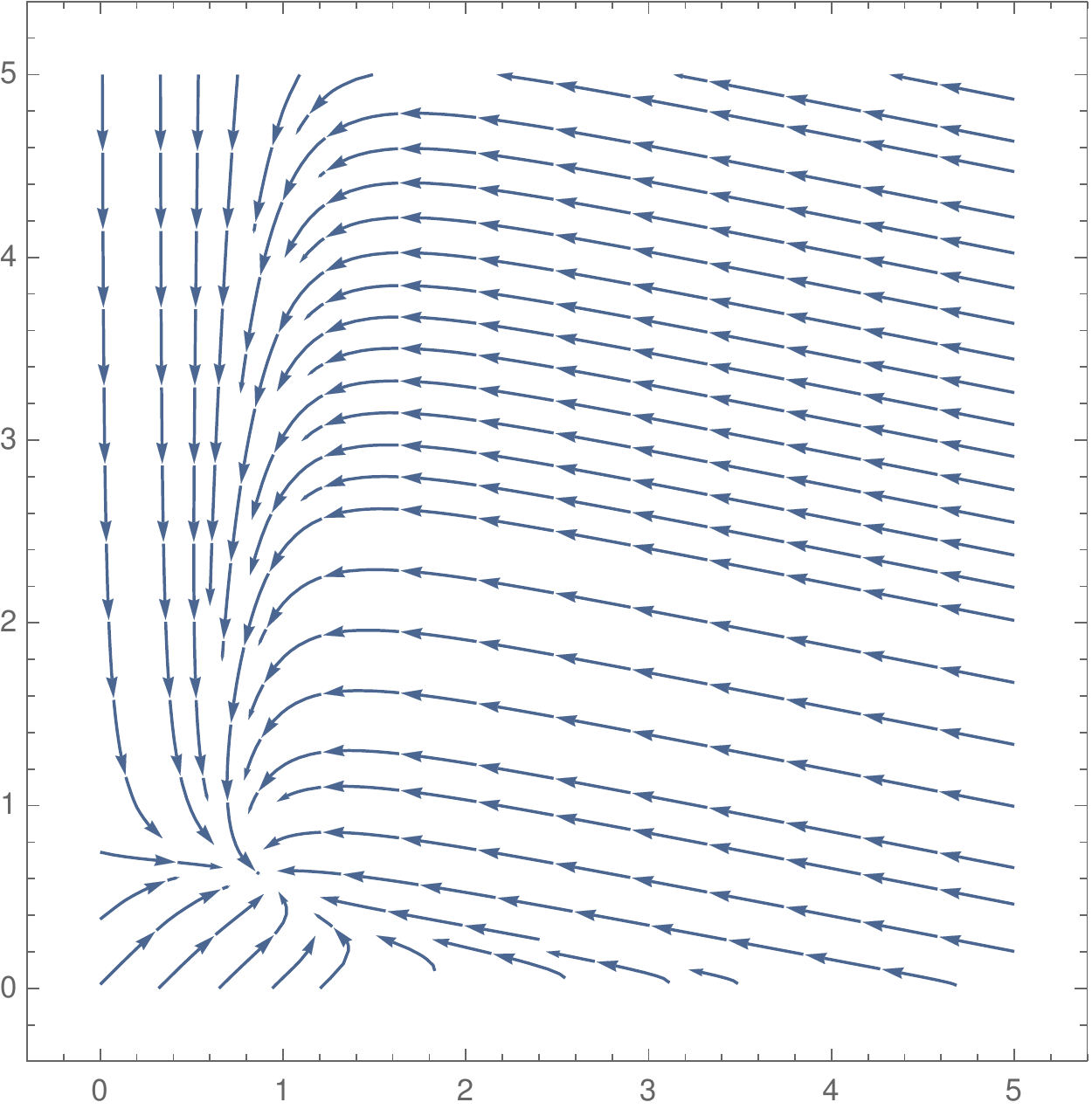}
  \caption{}
 \end{subfigure}
 \begin{subfigure}{.32\textwidth}
  \centering
  \includegraphics[width=.8\linewidth]{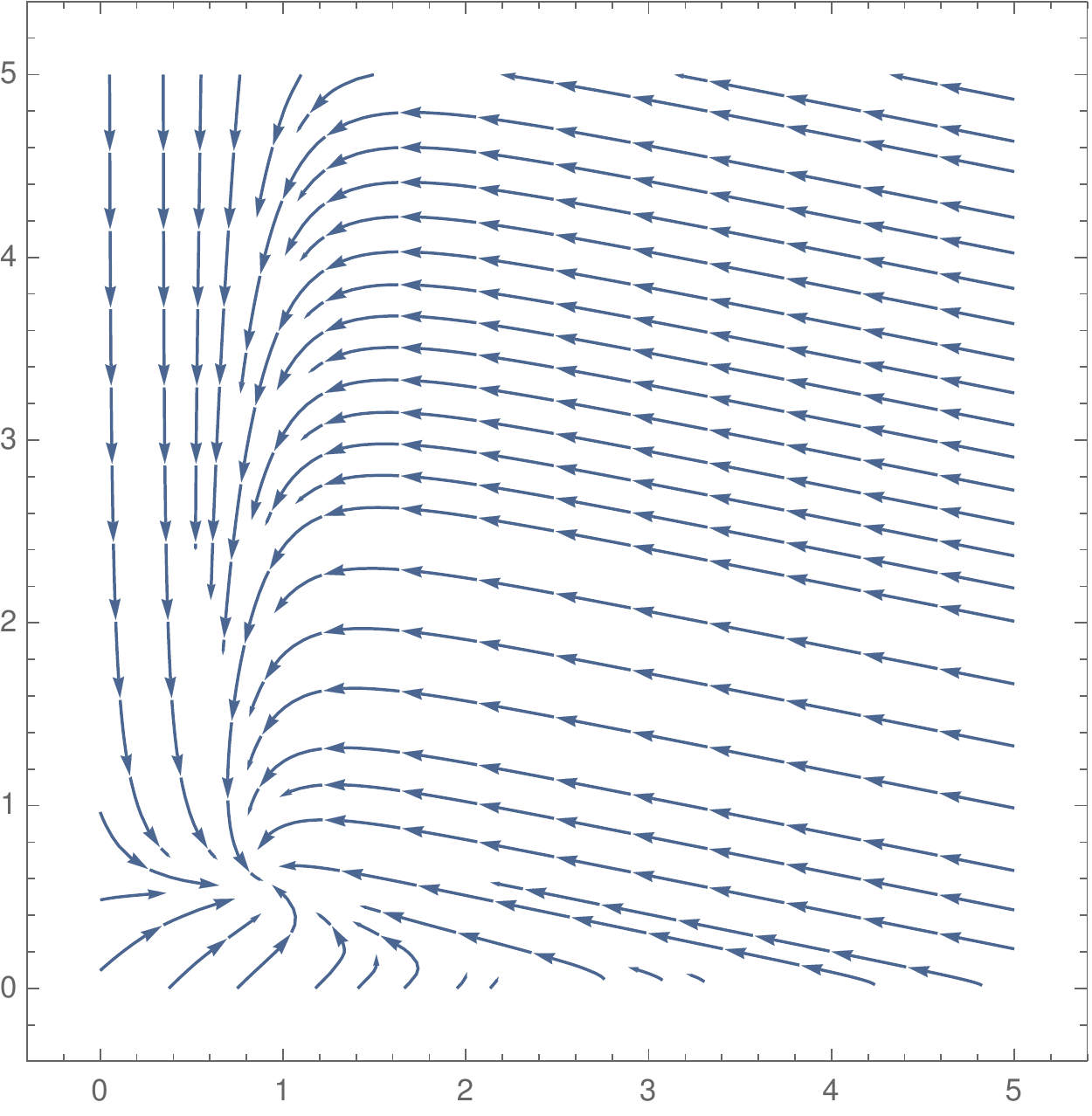}
  \caption{}
 \end{subfigure}
 \begin{subfigure}{.32\textwidth}
  \centering
  \includegraphics[width=.8\linewidth]{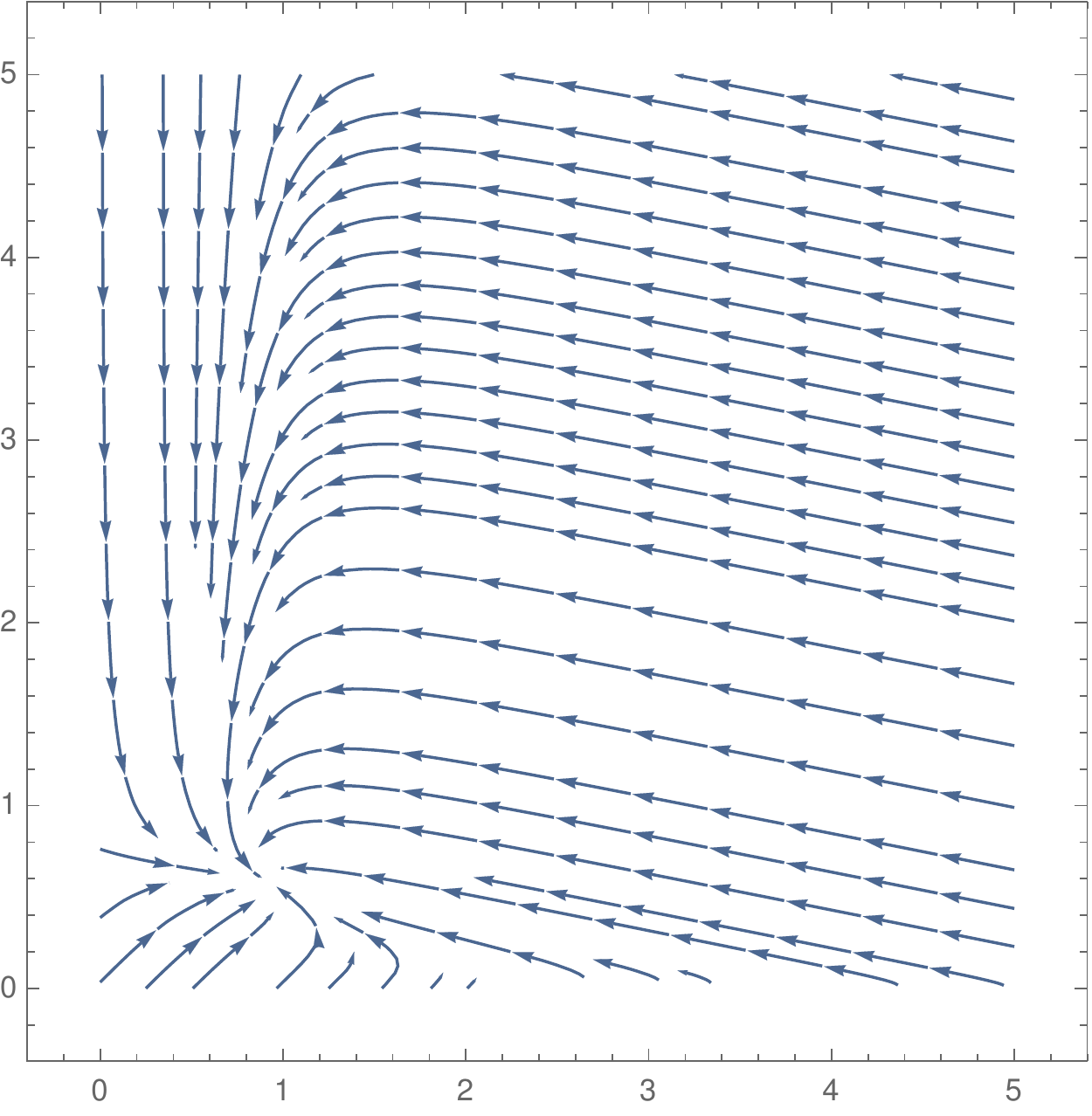}
  \caption{}
 \end{subfigure}
 \caption{Stream lines of the vector fields in \eref{e:odes} for the networks \abbr{crn}0 (a), \abbr{crn}1 (b) and \abbr{crn}2 (c). The vector fields are very similar, with the only noticeable differences close to the horizontal ($x_1$) axis in the three figures. Asymptotically, these differences vanish.}
 \label{f:f3}
\end{figure}

\paragraph{Behavior in the fluid limit regime} \label{s:fluidlimit}

Under the law of mass action, as the number of molecules in \eref{e:e0}--\eref{e:e2} goes to infinity the dynamics of their appropriately rescaled density obeys the system of ordinary differential equations \eref{e:ma}, \ie,
\begin{equ}\label{e:odes}
 \frac{\dd \ }{\dt} \begin{pmatrix} x_1\\x_2 \end{pmatrix} =  \begin{pmatrix} 1\\1 \end{pmatrix} + x_1^5x_2^2\begin{pmatrix}-5\\1 \end{pmatrix} + (x_2^3 + x_1^{\#} x_2) \begin{pmatrix}0\\-1 \end{pmatrix}\,,
\end{equ}
where $\#$ corresponds to the number of the \abbr{crn} under consideration (\emph{e.g.} $\# = 0$ for \abbr{crn}0) and without loss of generality we have assumed that $\kappa_r = 1$ for all $r \in \R$. The latter assumption will continue to hold throughout the paper.

The \abbr{crn}s defined above belong to the class of Asiphonic Strongly Endotactic networks, introduced in \cite{agazzi17,agazzi18}. This class of \abbr{crn}s is defined exclusively on structural properties of the networks. These properties on one hand that no subset of $\partial \RR^d$ is invariant under the flow of \eref{e:odes} (\ie that such networks can recover from the extinction of any of their species) and on the other hand the
asymptotic stability of \eref{e:odes} as $t \to \infty$ \cite{gopal13}. Furthermore, deviations of the appropriately rescaled stochastic systems \eref{e:e0}--\eref{e:e2} satisfy a large deviations principle in finite time \cite{agazzi17}.
Both these results are obtained by showing that at large concentrations the reactions dominating the sum \eref{e:generator} (the ones with input complexes on the vertices of the shaded polytope in \fref{f:f1}) contribute to pushing the state of the system towards  a compact set in phase space. We note that for \abbr{crn0}, \abbr{crn1} and \abbr{crn2} these dominating reaction coincide. Asymptotically, these networks display therefore the same behavior, as can be seen from the flow lines in \fref{f:f3}.

\paragraph{Behavior at finite size}

The probabilistic nature of the process $X_t$ requires a different definition of stability than the one given in the deterministic framework:
\dfn{\label{d:stability} For a compact set $\KK \in \Nn_0^d$, we denote the return time of the process $X_t$ to $\KK$ by $\tau_\KK := \inf\{t > 0~:~X_t \in \KK\}$ and say that a the process $X_t$ is \emph{positive recurrent} if $\Ex{x}{\tau_\KK} < \infty$, \emph{null recurrent} if $\px{x}{\tau_\KK < \infty} = 1$ but $\Ex{x}{\tau_\KK} = \infty$, and \emph{transient}  if $\px{x}{\tau_\KK < \infty} < 1$.}
To infer the stability properties of $X_t$ as defined in \dref{d:stability} for the systems \eref{e:e0}--\eref{e:e2}, we study similarly to \cite{hm11} whether the following Lyapunov stability condition is satisfied. Here $\indicator_A$ denotes the indicator function on the set $A$
\begin{cnd}[Stability]\label{cnd:stability}
 There exists a Lyapunov function $V$ with rate $\vphi$, \ie a continuous function $V\colon\mathbb N_0^d \to \RR$ with precompact sublevel sets such that
 \begin{equ}
  \LL V \leq K\,\indicator_\KK - \vphi(V)~,
  \label{e:phi}
 \end{equ}
 for some constant $K > 0$, a compact set $\KK$ and some monotone function $\vphi(x)\colon\RR\to \RR$~.
\end{cnd}
This condition guarantees that the process returns, on average, to a compact set in $\mathbb N_0^d$ by bounding form above its average speed towards that set.

Provided that the stability of the process $X_t$ has been established, we proceed to study its invariant measure, \ie a measure $\mu$ on $\Nn_0^d$ with $\mu\PP_t = \mu$\,. In particular, we are interested in the density of such measure and in the convergence under the flow defined by $\PP_t$ of another measure $\nu$ to it. Similarly to \cite{hm11} we obtain such result by combining \cndref{cnd:stability} with the following mixing condition:
\begin{cnd}[Mixing] \label{cnd:mixing} The level sets of $V$ are ``small'' enough, \ie for every $C_M>0$ and every $(x,y) \in \mathbb{N}_0^d \times\mathbb{N}_0^d$ such that $V(x) + V(y) \leq C_M$, there exists $\alpha > 0$ and $T>0$ such that
 \begin{equ}   \label{e:mixingcond}
  \tvn{\PP_T(x,\,\cdot\,) - \PP_T(y,\,\cdot\,)} \leq 2(1-\alpha)~.
 \end{equ}
\end{cnd}
Combined with \cndref{cnd:stability}, this Doeblin-like condition ensures the existence of a ``small'' attracting region of phase space where two processes with the same generator mix fast enough before fluctuating out of such region.

We are now ready to describe how, despite their similar asymptotic behaviour in the fluid limit regime, the three networks in \eref{e:e0}--\eref{e:e2} display significantly different behavior when modeled probabilistically. This is summarized in the following result

\begin{pro} \label{p:thm} The stochastic process $X_t$ for the network
 \begin{enumerate}[(a)]
  \item  \abbr{crn0} is positive recurrent and has a unique invariant probability measure $\mu_0$ satisfying $\int_{\Nn_0^2} \vphi(V(x)) \mu_0(\dd x) < \infty$. Any two initial point measures converge exponentially fast to one another, \ie for any $x, y \in \Nn_0^d$ there exists $\rho < 1$ and a positive constant $C$ such that for all $t > 0$ one has
  $$\tvn{\PP_t(x, \,\cdot\,) - \PP_t(y,\,\cdot\,)} \leq C \rho^t\tvn{\PP_0(x, \,\cdot\,) - \PP_0(y, \,\cdot\,)}\,.$$
  Consequently, we have exponential convergence to the invariant measure, \ie for any $x \in \Nn_0^d$ there exists $\rho' < 1$ and a positive constant $C'$ such that for all $t > 0$ one has
  $$\tvn{\PP_t(x, \,\cdot\,) - \mu_0} \leq C' \rho'^t\tvn{\PP_0(x, \,\cdot\,) - \mu_0}\,.$$
  \item  \abbr{crn1} is null recurrent and has a unique (up to scalar multiplication) $\sigma$-finite invariant measure $\mu_1$ satisfying $\int_{\Nn_0^2} \vphi(V(x)) \mu_1(\dd x) < \infty$. \cf{Furthermore, for any $x,y \in \Nn_0^d$ there exists $C'' > 0$ such that
  $$\tvn{\PP_t(x, \,\cdot\,) - \PP_t(y, \,\cdot\,)} \leq \frac{C''}{H_\vphi^{-1}(t)}\,,$$
  where $H_\vphi^{-1}(t)$ is the inverse function of $H_\vphi(u) := \int_1^u \vphi(s)^{-1}\dd s$\,.}
  \item  \abbr{crn2} is transient.
\end{enumerate}
\end{pro}

We establish the above result in two steps. The first is to construct a Lyapunov function $V$ satisfying \cndref{cnd:stability} when $X_t$ is not transient, while the second is to combine \cndref{cnd:stability} and \cndref{cnd:mixing} to obtain the desired estimates on $\mu$.


 The first step is executed by studying the dominant behaviour of the generator $\LL$ for large values of $x$. Indeed, as shown in \cite{agazzi18}, the phase space can be divided radially (in a certain set of coordinates) into dominance regions where $\LL$ acquires asymptotically a particularly simple form. We use this simplification to construct a Lyapunov function $V$ satisfying \eref{e:phi} in each of these regions separately. Finally, we handle the gluing of the locally defined candidate Lyapunov function. This is done by showing that, under certain convexity conditions on the glued Lyapunov function at the interface between two adjacent regions, \eref{e:phi} is automatically satisfied at that interface.\\
 The strategy adopted for the local construction of $V$ to simplify the verification of \eref{e:phi} in a given dominance region $\TT \subset \RR^2$ pivots on the monomial form of the asymptotic rates $\lambda_r$. Indeed, such rates scale homogeneously under any scaling transformation. This is true in particular for transformations leaving subsets of $\TT$ invariant. Hence, constructing $V$ to also scale homogeneously under such transformations allows to exclude most of the summands in $\LL V$ by power counting, reducing the right hand side in \eref{e:phi} to a monomial in the scaling variable, whose evaluation is immediate. Taking also $h$ to scale homogeneously under the same transformation reduces the left hand side of \eref{e:phi} to the same form, significantly simplifying the verification of the desired inequality.

 In the second step, the results of \pref{p:thm} are obtained through \cndref{cnd:stability} and \cndref{cnd:mixing} by similar arguments to the ones developed in \cite{hm11}.

The paper is structured as follows: in Section~\ref{s:boundary} we study the behavior of the process $X_t$ for small values of $(X_t)_2$. This is the region where the dynamics of the three examples \eref{e:e0}-\eref{e:e2} are radically different from one another. This allows to prove \pref{p:thm} (c). In Section~\ref{s:Lyapunov} we turn to the piecewise construction of Lyapunov functions for the study of the stability of the two remaining examples. In Section~\ref{s:Lyapunov-2} we study the behavior of $\LL$, which we use in Sections~\ref{s:Lyapunov-1}--\ref{s:Lyapunov1} to construct $V$ satisfying \eref{e:phi} as described above. Section~\ref{s:Lyapunov2} is devoted to the patching of local Lyapunov functions. Finally, in Section~\ref{s:convergence} we prove the results about existence of and convergence to the invariant measure.

\section{Stability of CRNs} \label{s:boundary}

As explained in Section~\ref{s:fluidlimit} and more precisely in \cite{agazzi17}, in the limit of large number of molecules the dynamics of \eref{e:e0}--\eref{e:e2} is stable under the appropriate scaling. This result relies on the approximation of the reaction rates \eref{e:rates} with monomials by the Stirling formula. This approximation breaks down when at least one component of $X_t$ is $\OO(1)$, \ie $X_t$ is close to the boundary of $\RR^d$. In this regime the stability results above are no longer valid in general and more careful analysis is needed to infer the dynamical behavior of the networks at hand.

\subsection{Boundary Dynamics}

In this section we consider the behavior of the three \abbr{crn}s when $(X_t)_2 = \OO(1)$. Recall that, by definition, some of the discrete rates $\Lambda_r(x)$ of the process $X_t$ will vanish in a neighborhood of the boundary $\{x_2 = 0\}$. More precisely, defining throughout $\TT_{01}(n) := \{(x_1,x_2) \in \mathbb N_0^2~:~x_2 < n\}$ we have that the set of reactions  $\{r_2,r_3,r_4\}$ will have vanishing rate in $\TT_{01}(1)$, while reactions in $\{r_3,r_4\}$ will do the same in $\TT_{01}(2)$. This intermittent behavior implies that within $\TT_{01}(2)$ noise dominates the dynamics of the process $X_t$ and its effects on the three networks differ significantly, as we see below. We use this to prove part (c) of \pref{p:thm}, and to characterize the exit distribution from this noise-dominated region for \abbr{crn0} and \abbr{crn1}.

\subsubsection{\abbr{crn2}}

Choosing without loss of generality the rate constants $\kappa_1 = \kappa_2 = 1$\,, the reaction rates of the network for $x=(x_1,x_2) \in \TT_{01}(2)$ are
\begin{equ}
 \Lambda_1(x) = 1\,, \quad \Lambda_2(x) = \indicator_{\{x_2 \geq 1\}} x_1^2 + \OO(x_1)\,, \quad \Lambda_3(x) = 0\,,\quad \Lambda_4(x) = 0~.\label{e:rates2}\end{equ}
 Throughout, we define the Markov chain $X_n$ associated to the process $X_t$ stopped when exiting $\TT_{01}(2)$ as 
 \begin{equ}
   \label{e:makovchain}X_n :=  X(\tau_n) \qquad \text{where } \qquad \tau_n = \inf\{t > \tau_{n-1} ~:~X(t) \neq X(\tau_{n-1}) \text{ or } X(t) \not \in \TT_{01}(2) \}\,.
\end{equ}
($\tau_n$ is finite with probability one since the jump rate is positive at any point in $\TT_{01}(2)$).
 The possible transitions of $X_n$ are:
 \begin{myitem}
  \item for $X_n = (x_1, 0)$ we jump to $X_{n+1} = (x_1 + 1,1)$  with probability $1$~,
  \item for $X_n = (x_1, 1)$ we jump up ($\uparrow$) to $X_{n+1} = (x_1 + 1,2)$  with probability $p_\uparrow(x_1) = 1/(x_1^2 + 1)$, or down ($\downarrow$)  to $X_{n+1} = (x_1,0)$  with probability $p_\downarrow(x_1) = 1- p_\uparrow(x_1) = x_1^2/(x_1^2 + 1)$~,
  \item for $X_n = (x_1, 2)$ the process is stopped, \ie $X_{n+1} = (x_1,2)$  with probability $1$~.
\end{myitem}

\noindent We are now ready to prove \pref{p:thm} (c).

\begin{cla}\label{cla:transiant}
 The Markov process $X_t$ associated to the \abbr{crn} \eref{e:e2} is transient.\label{c:div2}
\end{cla}
\begin{proof}
 By irreducibility of $X_t$ we can choose without loss of generality $X_0 = (k_0,0)$ for $k_0 > 0$. As $X_t$ and $X_n$ are equivalent up to time change,
it is sufficient to show that the associated Markov chain $X_n$ is transient.
 To do so, we bound from below the probability that $X_n$ ``oscillates forward'' never leaving the tube $\TT_{01}(2)$: For $k \in \Nn$ we denote by $\uparrow_k$ and $\downarrow_k$ the event of jumping up and down, respectively, from site $(k,1)$. Observe that for any $n > k_0>0$ we have
 $$\px{(k_0,0)}{X_n \in \TT_{01}(2)} =  \p{\bigcap_{k = k_0+1}^n \downarrow_{k} } \geq 1-\p{\bigcup_{k = k_0+1}^n \uparrow_k} = 1- \sum_{k = k_0}^n p_\downarrow(k) = 1 - \sum_{k = k_0}^n \frac{1}{k^2 + 1}~.$$
 Hence, the claim follows by taking the limit $n \to \infty$ in the display above and bounding the right hand side from below and away from $0$ upon choosing $k_0$ large enough.
\end{proof}

\subsubsection{\abbr{crn0} and \abbr{crn1}}

We proceed to characterize the exit distribution of $X_t$ from $\TT_{01}(2)$ for \abbr{crn0} and \abbr{crn1}. For these networks, the reaction rates are identical to those in \eref{e:rates2} except for $\Lambda_2(\,\cdot\,)$, which reads,
\begin{equ}
 \text{\abbr{crn0}: } \quad \Lambda_2(x) = 1\,,\qquad \text{\abbr{crn1}: } \quad \Lambda_2(x) = \indicator_{\{x_2 \geq 1\}} x_1\,.
\end{equ}
Consequently, the jumping probabilities of the associated Markov chain $X_n$ remain unchanged except when $x_2 = 1$. Assuming without loss of generality that $\kappa_r = 1$ for all $r \in \R$, we have for $X_n = (x_1,1)$  that
\begin{equ} \label{e:probasjumps}
  \text{\abbr{crn0}: } \quad p_\uparrow(x_1) = \frac{1}{2}\,, \quad p_\downarrow(x_1) = \frac{1}{2}\,,\qquad \text{\abbr{crn1}: } \quad p_\uparrow(x_1) = \frac{1}{x_1 + 1}\,, \quad p_\downarrow(x_1) = \frac{x_1}{x_1 + 1}\,.
\end{equ}
Therefore, the exit distribution of $X_t$ from $\TT_{01}(2)$ is the same as for $X_n$. Furthermore  for $X_0 = (k_0,0)$ and $b \geq k_0$,  we have
\begin{equ}\label{e:exitdistro}
  \px{(k_0,0)}{X_\infty = (b+1,2)} := \lim_{n \to \infty}\px{(k_0,0)}{X_n = (b+1,2)} =  \p{ \uparrow_{b} \cap \bigcap_{k = k_0+1}^{b-1} \downarrow_k } = p_\uparrow(b) \prod_{k = k_0+1}^{b-1} p_\downarrow(k)\,.
\end{equ}
Combining \eref{e:probasjumps} and \eref{e:exitdistro} we obtain for \abbr{crn0}
\begin{equ}
 \px{(k_0,0)}{X_\infty = (b+1,2)} = (1-a) a^{b-k_0}\,,
\end{equ}
where we have defined $a = \p{\downarrow} = 1/2\,$ the probability of jumping down on level $x_2 = 1$. Similarly for \abbr{crn1} we have 
\begin{equ}\label{e:pupcrn1}
 \lim_{n \to \infty}\px{(k_0,0)}{X_\infty = (b+1,2)} = \frac{1}{b+1} \prod_{k = k_0}^{b-1} \frac{k}{k+1} = \frac{k_0}{b(b+1)}\,.
\end{equ}
\rmk{\label{r:nullrecurrent} Using \eref{e:pupcrn1} and defining $k_0 := (x_0)_1$ we bound from below the expected hitting time $\tau_\KK$ of a compact set $\KK \subseteq \RR^d$ by
\begin{equs}
  \Ex{x_0}{\tau_\KK} \geq \Ex{x_0}{\inf\{\tau~:~(X_\tau)_{2}>1\}} \geq \sum_{k = k_0}^{\infty}\px{x_0}{X_\infty = (k+1,2)}\sum_{j = k_0}^{k} \E{\Delta\tau_{(j,0)}}
  \geq  x_0\int_{k_0+1}^\infty \frac{k-x_0}{k (k+1)} \dd k\,,
\end{equs}
where we used that the jumping time $\delta \tau_{x}$ at $(j,0)$ has constant expectation. Hence, for \abbr{crn1} we have $\E{\tau_\KK} = \infty$. In light of this, in order to prove null recurrence it remains to show that $\px{x}{\tau_\KK < \infty} = 1$. By \cite[Theorem 12.3.3]{meyn12} this is the case if there exists a Lyapunov function satisfying \eref{e:phi}. We construct such a function in the upcoming section.}

\section{Local Lyapunov functions}
\label{s:Lyapunov}

In this section we study the stability of \abbr{crn0} and \abbr{crn1} by introducing a family of Lyapunov functions $V$ that verifies \cndref{cnd:stability}. We outline below the intuition behind the construction of $V$ and will carry out such construction in full detail in the subsequent sections.

To establish \eref{e:phi} we notice that, by boundedness of the function $V$ and of the rates \eref{e:rates} on compact sets, it is sufficient to have
\begin{equ}\label{e:phi3}
  \LL V(x) \leq - \vphi \circ V(x) \qquad \forall\, x \in \Nn_0^d \setminus \KK\,,
\end{equ}
for a compact set $\KK$ large enough. In other words we are interested in the asymptotic behavior of the left hand side of \eref{e:phi3}. To explore it we introduce the following family of scaling transformations
\begin{df}[Scaling transformations]\label{d:trafo}
For any vector $w = (w_1,w_2) \in S^1$ and scaling parameter $l \in \Rr_+$ we define the family of transformations $\slw$  by
\begin{equ}\label{e:trafo}
  \slw(x) := (l^{w_1}x_1,l^{w_2}x_2)\quad\text{for}\quad x=(x_1,x_2) \in \Nn_0^d .
\end{equ}
Furthermore, we say that a function $\phi\colon\RR^2 \to \Rr$ scales homogeneously under $\slw$ if $\phi \circ \slw (\,\cdot\,) = l^{\delta} \phi(\,\cdot\,)$ for some $\delta \in \Rr$. In this case we write $\phi \ssim{w} l^\delta$. \label{d:homogeneous}
\end{df}

\noindent For any $w \in S^1$ we then obtain \eref{e:phi3} in a region $\TT$ that is invariant under \eref{e:trafo} by showing that for $l$ large enough
\begin{equ} \label{e:SLW!}
  \sum_{r \in \R} \Lambda_r(\slw (x)) \Delta_r V(\slw (x)) \leq - \vphi \circ V(\slw(x))\,.
\end{equ}
for all $x \in \RR^d$ with $\|x\|_1 \leq 1$ and $\slw (x) \in \Nn_0^d$, where we define throughout the difference operator $\Delta_r f(x) := f\pc{x + c_r} - f(x)$.

Let throughout $S_+^1$ be the open positive orthant of $S^1$. Fixing $w \in S_+^1$ recall that the rates \eref{e:rates} scale, to leading order, like the ones in \eref{e:ma}:
\begin{equ} \label{e:approxrates}
  \Lambda_r (\slw(x)) = \kappa_r \prod_{i = 1}^d \binom{l^{w_i}x_i}{(c_\inn^r)_i}(c_\inn^r)_i! = \kappa_r l^{\dtp{w, c_\inn^r}} \prod_{i = 1}^dx_i^{(c_\inn^r)_i} +  \OO(l^{\dtp{w, c_\inn^r}}) = \lambda_r(\slw (x)) + \OO(l^{\dtp{w, c_\inn^r}})\,,
\end{equ}
and that all the rates $\lambda_r$ scale homogeneously under $\slw$. Under the key assumption that the same holds for $V$, \ie provided that there exists a function $\obar V \in C^1(\RR^2)$ and $\delta\colon S^1 \to \RR$ such that
\begin{equ}\label{e:homogeneous}
  V(\slw (x)) = l^{\delta(w)} \obar V(x)\,,
\end{equ}
we show that for all $r \in \R$ we can write
\begin{equ} \label{e:homogeneousop}
  \Delta_r V(\slw (x)) = l^{v_r(w)} \obar \Delta_r \obar V(x) + \OO(l^{v_r(w)})\,,
\end{equ}
 for an operator $\obar \Delta_r\colon C^1(\RR^d) \to C^0(\RR^d)$\, and a function $v_r\colon S^1 \to \Rr$\,.
This allows to obtain the desired result by a scaling argument: writing the left hand side of \eref{e:SLW!} to leading order in each of its summands as
\begin{equ} \label{e:scaling}
 \LL V(\slw(x)) = \sum_{r \in \R} \pc{l^{\dtp{w, c_\inn^r} + v_r(w)}\lambda_r(x) \obar \Delta_r \obar V(x) + \OO(l^{\dtp{w, c_\inn^r}})} =  l^{\delta'(w)} h(x) + \OO(l^{\delta'(w)}) \,,
\end{equ}
where
\begin{equ}\label{e:deltass}
\delta'(w) := \max\{\delta_r'(w)~:~r \in \R\} \qquad \text{for} \qquad \delta_r'(w) := \dtp{w,c_r^\inn} + v_r(w)
\end{equ}
and
\begin{equs} \label{e:h}
  h(x) :=  T_i \obar V(x)\,,\qquad  \text{for} \qquad T_i f (x) := \sum_{r \in \R_w}\lambda_r(x) \obar \Delta_r f(x)\,,
\end{equs}
and $\R_w :=  \{r \in \R\colon  \delta_r'(w) \geq \delta_{r'}'(w) \,\forall r' \in \R\}$, we immediately obtain \eref{e:phi3} by combining \eref{e:homogeneous} with \eref{e:scaling} and by defining
\begin{equ}\label{e:phiscale}
  \vphi(x) := \Ch x^{\delta'(w)/\delta(w)}\, \qquad \text{for}\quad  \Ch \in (0,1)\,.
\end{equ}
\begin{remark}
  Assumption \eref{e:homogeneous} emerges naturally from the structure of the problem at hand, as encoded by the generator \eref{e:generator}. Indeed, by the scaling of the rates \eref{e:rates} as monomials for all $w \in S^1$, this assumption allows to consider the left hand side of \eref{e:phi3} as a polynomial in $\slw$. Consequently, for any $w \in S^1$ the study of $\LL V$ reduces to identifying the term dominating the polynomial under that scaling. This establishes a connection with the domain of tropical geometry \cite{strumfels15}.
\end{remark}
In the following sections we realize the program outlined above. First of all we characterize the operator $\obar \Delta$ for which \eref{e:homogeneous} holds, and divide the space into regions $\{\TT\}$ where the leading term in the approximation of $\LL V$ is constant. We then locally construct a region-specific Lyapunov pair $(\obar V , h)$ satisfying the definition \eref{e:scaling} by solving the Poisson equation \eref{e:h} in each region, \ie
\begin{equ}
 \begin{cases}
  T_i \obar V_i(x) = -h_i(x)   & \text{for } x \in \TT_i         \\
  \obar V_i(x) = \obar V_j(x) & \text{for } x \in \partial \TT_i
\end{cases}\,.\label{e:bvp}
\end{equ}
for a function $h \ssim{w} l^{\delta'(w)}$ and boundary conditions $\obar V_j \ssim{w} l^{\delta(w)}$. This way we enforce \eref{e:scaling} to the first order in the scaling parameter $l$. By the leading order expansion we expect such candidate Lyapunov functions to solve \eref{e:phi3} with $\vphi$ as in \eref{e:phiscale}. We verify that this condition is indeed satisfied in the last paragraph of this section.

\subsection{Scaling of the generator} \label{s:Lyapunov-2}


Assuming that {for $w \in \obar{S_+^1}$} the function $V(x)$ satisfies \eref{e:homogeneous} we have that
\begin{equs} \label{e:SLW}
  \LL V(\slw(x)) &= \sum_{r \in \R} \Lambda_r(\slw (x)) \pc{V(\slw (x) + c^r) - V(\slw (x)}\\ &= \sum_{r \in \R} l^{\dtp {c_\inn^r ,w}+v_r} \lambda_r(x) \pc{ \obar V(x + (\slw)^{-1}( c^r)) - \obar V(x)} + \OO(l^{\dtp{c_\inn^r, w}+v_r})\,.
\end{equs}
Now, for $w \in S_+^1$ we can approximate to leading order the difference terms in \eref{e:SLW} by partial derivatives. Indeed, for each $r \in \R$ we expand $\obar V$ as
\begin{equ}\label{e:taylor}
  \obar V\pc{x + (\slw)^{-1}( c^r)}  = \obar V\pc{x_i +  l^{-w_{i}}c^r_{i}} = \obar V(x) + \sum_{i \in I_r} l^{-w_i}{c^r_{i} \partial_i \obar V(x)} + \OO(l^{-w_i})\,.
\end{equ}
where $\partial_i$ denotes a partial derivative in direction $i \in \S$ and we have defined the index set $I_r := \min\{i \in \supp c^r~:~ w_i \leq w_j \forall j \neq i\}$. We make this statement precise in \lref{l:taylor2} below.
We see that the dominant term in \eref{e:SLW!} depends on the chosen vector $w$. Therefore we can divide $S^1_+$ into regions characterized by different dominating terms in \eref{e:SLW!}. In this case we write:
\begin{equ}
 \begin{cases}
  \LL \sim T_1 :=   x_1^5x_2^2   (-5 \partial_1)               & \text{for } w \in \WW_1 := \{w \in S^1~:~w_2 > w_1 > 0\} \\
  \LL  \sim T_2 :=  x_1^5x_2^2 \pc{   \partial_2 -5 \partial_1 } & \text{for } w \in \WW_2 := \{w \in S^1~:w_2 = w_1 = 1/\sqrt 2\}  \\
  \LL  \sim T_3 := x_1^5x_2^2  \partial_2                      & \text{for } w \in \WW_3 := \{w \in S^1~:0 < w_2 < w_1\}
\end{cases}\,. \label{e:transportgen}
\end{equ}
\begin{remark}\label{r:asgenerator}
  We see that the dominant terms in \eref{e:transportgen} and the dominance regions are the same as the ones of the deterministic transport generator $T f(x) := \sum_{r \in \R} \lambda_r(x) \dtp{c^r, \nabla f(x)}$, which for \abbr{crn0} reads
\begin{equ}
 T := \pc{l^{-w_1}\partial_1 + l^{-w_2}\partial_2} + l^{-w_2} (-\partial_2) + l^{5w_1 + 2w_2} x_1^5x_2^2 \pc{l^{-w_1}  (-5 \partial_1) + l^{-w_2} \partial_2} + l^{2w_2} x_2^3 (- \partial_2)~.\label{e:Ltransport}
\end{equ}
In this sense, for $w \in S_+^1$ the discrete generator \eref{e:generator} is well approximated under $\slw$  by its deterministic, continuous counterpart $T$\,. For $w \in \{(1,0), (0,1)\}$ this approximation is not possible because the jumping nature of the process becomes dominant, as we will see in the upcoming sections.
\end{remark}

\begin{figure}
 \centering
\begin{subfigure}{.45\textwidth}
  \centering
  \def\svgwidth{1\textwidth}
  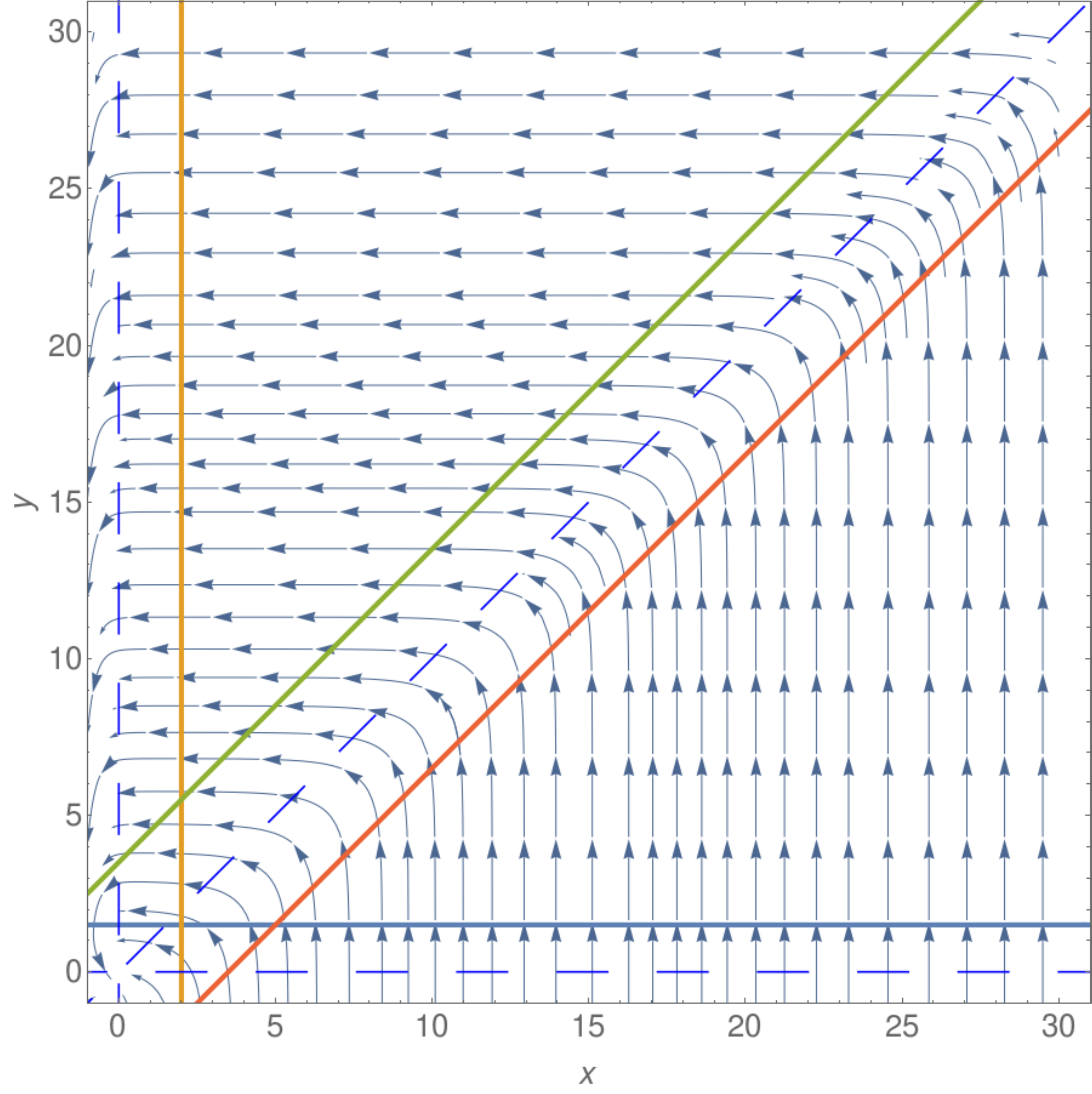
  \caption{}
\end{subfigure}
\qquad
\begin{subfigure}{.45\textwidth}
 \centering
 \def\svgwidth{1\textwidth}
 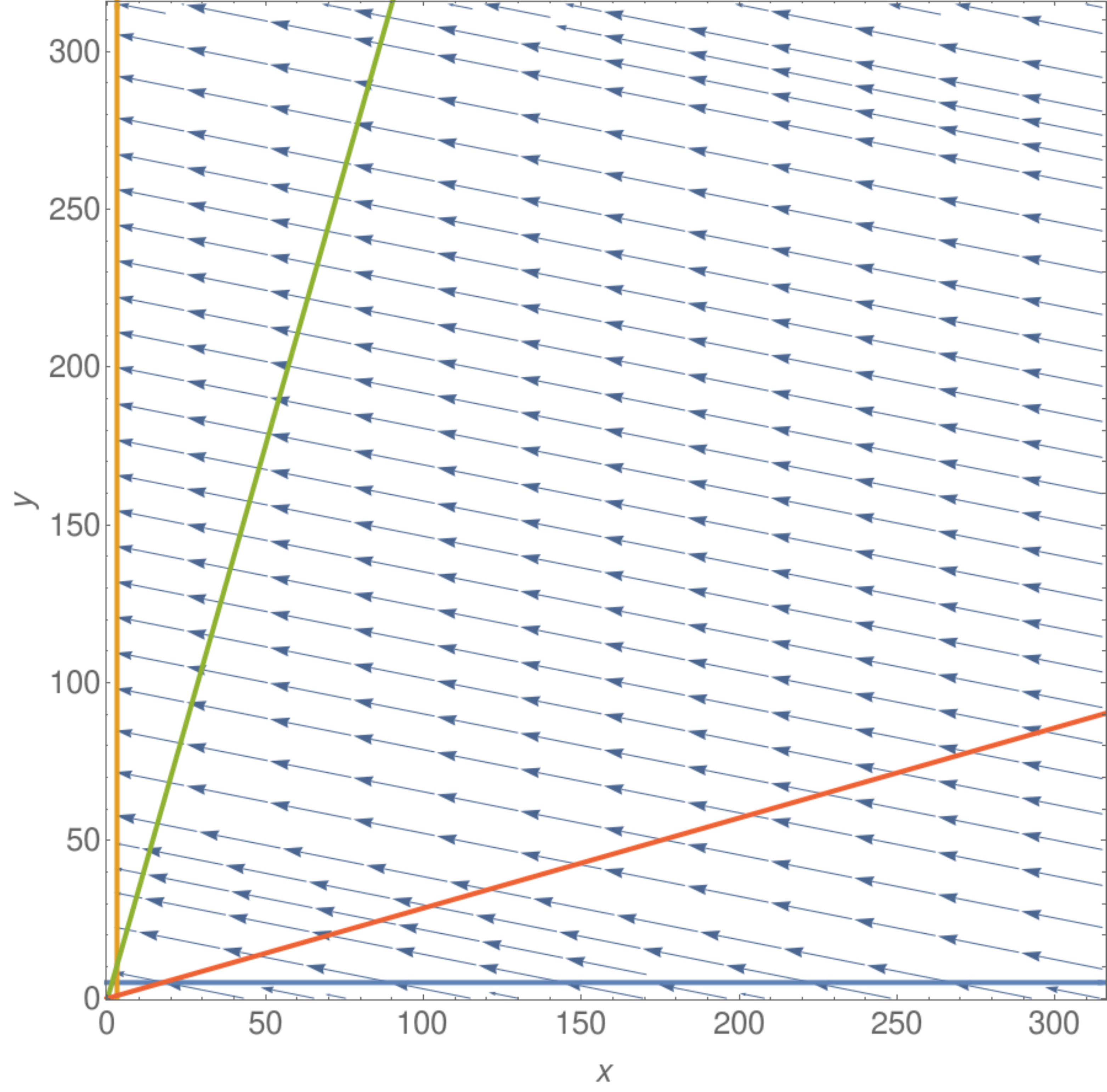
 \caption{}
\end{subfigure}
 \caption{The vector field for \eref{e:Ltransport} under the $\slw$ scaling for large $l$ in the coordinates $(l,w)$ for $w \in S_+^1$ (a) and in the original coordinates $(x_1,x_2)$ (b). In (a) dashed lines denote the radial directions associated to $\WW_0, \WW_2, \WW_4$. The corresponding neighborhoods $\WW_0(2), \WW_2(4), \WW_4(2)$ are mapped by the component-wise exponential into the transition regions $\TT_0, \TT_2, \TT_4$ respectively. Solid lines denote transitions between different dominance regions in both coordinate systems: in (b) $\TT_{01}$ in blue, $\TT_{12}$ in red, $\TT_{23}$ in green and $\TT_{34}$ in yellow). }
 \label{f:change}
\end{figure}

  It is apparent in \eref{e:transportgen} that the point sets $\WW_0 := \{(0,1)\}, \WW_2 = \{1/\sqrt 2 (1,1)\}, \WW_4 := \{(1,0)\}$ limit intervals of $w$ where one single term dominates the generator. These points uniquely identify radial lines in $\Rr^2$ through the polar coordinate system with parameters $(l,w) \in (\RR, S^1)$. In the neighborhoods of these lines, defined throughout as
  \begin{equ}\label{e:wneighborhood}
    \WW_{i}(\xi) := \pg{(\theta, w) \in \Rr_{\geq 1} \times S^{d-1}~:~\log(\theta) \inf_{w' \in \WW_i} \| (w-{w'})\|_2 \leq \xi}\,,
  \end{equ}
   we observe a transition between dominant terms of the generator, as depicted in \fref{f:change}~(a). Such neighborhoods and their complement define a partition of $\Rr^2$. This partition can be mapped through the component-wise exponential function to a partition of $\RR^2$ into dominance regions of the generator, as displayed in \fref{f:change}~(b).
   \begin{remark}\label{r:agazzi18}
     The partition $\WW_{0}, \dots, \WW_{4}$ defined above corresponds to the partition $\WW_1^* \cup \WW_2^*$ introduced in \cite{agazzi18}. (See  \cite{agazzi18} for the definition of $\WW_1^*$ and $\WW_2^*$.)  This construction is therefore naturally generalizable to a higher-dimensional framework. Furthermore, the regions defined in \eref{e:wneighborhood} correspond, asymptotically in $l$, to the sets introduced in \cite[Definition~4.16]{agazzi18}, and the geometric results of \cite[Lemma 4.26]{agazzi18} therefore directly apply to the problem at hand.
   \end{remark}
We denote the regions in $\RR^2$ corresponding to the dominant behaviors in \eref{e:generator} by $\TT_i$ for $i \in (0, \dots, 4)$.
We further define, throughout, the boundaries separating two regions $i$ and $j$ by $\TT_{ij} := \obar \TT_i \cap \obar \TT_j$, where $\obar A$ denotes the closure of the set $A$. In our case these sets can be written as
\begin{equ}\label{e:boundaries}
 \TT_{01} := \{x_2 = b_0\}\,,\quad \qquad \TT_{12} :=\{x_2 =  x_1/b_1\} \,,\quad \qquad \TT_{23} :=\{x_2 =  b_1 x_1\}\,,\quad \qquad \TT_{34} :=\{x_1 = b_2\}\,,
\end{equ}
where the parameters $b_0, b_1, b_2 \in (1, \infty)$ will be fixed to our convenience at a later point of the analysis. The dominant behaviors of the generator and the above definitions are summarized in \fref{f:change}.

\begin{remark} The boundaries separating different dominating regions of the generator are so-called toric rays and they partition the phase space into a tropical fan \aa{cite}. A similar asymptotic partition was used in \cite{agazzi18, gopal13} always for the study of \abbr{crn}s, establishing a connection between this subject and the one of tropical geometry.\end{remark}

\subsection{Construction of the Lyapunov function} \label{s:Lyapunov-1}

We now use the partition developed in the previous section to construct a Lyapunov function satisfying \eref{e:h} and scaling homogeneously as required in \eref{e:homogeneous}. More precisely we construct such a Lyapunov function $\obar V_i(x)$ in each $\TT_i$ by solving the Poisson problem \eref{e:bvp} with $\partial \TT_i = \TT_{ij}$ for all $j$ with $\obar \TT_{i} \cap \obar \TT_{j} \neq \emptyset$.
The reason for this choice of Lyapunov function is twofold. On one hand, since both $T_i$ and our choice of $h_i$ and boundary conditions $\obar V_j$ scale homogeneously under $\slw$ for $w \in \WW_i$, we expect the candidate Lyapunov function to also scale homogeneously under the same set of transformations, thereby satisfying the key assumption \eref{e:homogeneous} of Section~\ref{s:Lyapunov-2}. On the other hand, as the function $\obar V_i$ satisfies \eref{e:h} by construction, by the dominance of $T_i$ in $\TT_i$ and the discussion of Section~\ref{s:Lyapunov-2} we expect it to satisfy \eref{e:phi3} with $\vphi$ as in \eref{e:phiscale}.

We solve \eref{e:bvp} by the method of characteristics, \ie by direct integration of the right hand side:
\begin{equ}
 \obar V_i(x) = \Ex{x}{\obar V_j(X_{\tau_i })} + \Ex{x}{\int_0^{\tau_i} h_i \pc{X_t} \dt}\,,\label{e:prop}
\end{equ}
where $\tau_i := \inf\{t~:~X_t \not \in \TT_i\}$ is the exit time from $\TT_i$ and expectations are taken with respect to the dynamics of $X_t$  with generator $T_i$ from \eref{e:bvp}. We proceed to calculate \eref{e:prop} by considering different regions of phase space separately. In each of these regions we will assume that $h$ scales homogeneously as in \dref{d:homogeneous} under a scaling that depends on the region where they are defined.

\subsubsection{Priming region: $\TT_4$}
\label{s:priming}

We start by studying the dynamics of $X_t$ in $\TT_4$. As displayed in \fref{f:change}(a) in this region of phase space the process $X_t$ is rapidly converging to a compact set. For this reason, we call $\TT_4$ the \emph{priming} region. We establish \eref{e:phi3} asymptotically by using the transformation from \dref{d:trafo} that leaves $\TT_4$ invariant, \ie $\slw$ with $w = (0,1)$.
We seek to obtain the Lyapunov function by solving \eref{e:bvp} with the assumption that both $h_4$ and the boundary condition $V_4^*$ on $\TT_4^* := := \{(x,y) \in \mathbb N_0^2~:~x = 0\}$ scale homogeneously under such transformation.
In order to realize this program we first identify the dominant terms of the generator $\LL$ under the scaling $\slw$. We proceed similarly to \eref{e:SLW} and assuming that \eref{e:homogeneous} holds we obtain
\begin{equ}\label{e:l0}
 \LL V(\slw x) = l^{3} x_2^3 \pc{\obar V(x_1,x_2-l^{-1})-\obar V(x)} +  l^2 \binom{x_1}{5} x_2^2\pc{\obar V(x_1-5,x_2)-\obar V(x)} + \OO(l^2).
\end{equ}
where we have not approximated the binomial coefficients involving $x_1$ as \eref{e:approxrates} does not hold in this regime and used that the first difference term on the right hand side scales as $l^{-1}$to leading order by Taylor expansion. In light of this we write the asymptotic generator in $\TT_4$ as
\begin{equ}\label{e:t4}
  T_4 f(x) = x_2^3 \pc{f(x_1,x_2-1)-f(x)} + \binom{x_1}{5} x_2^2\pc{f(x_1-5,x_2)-f(x)}\,.
\end{equ}
\begin{remark}\label{r:differencegen}
  Because $w_1 = 0$\,, in this region \eref{e:approxrates} and \eref{e:taylor}  cannot be used to approximate \eref{e:l0} by a first order differential operator with monomial rates as in \eref{e:SLW}. In particular, despite being significantly simplified with respect to its original form, the operator $T_4$ remains the one of a jump Markov process. In this sense, the discrete nature of the problem at hand is ``felt'' by the Lyapunov function only close to $\partial \RR^2$.
\end{remark}
We now proceed to propagate the Lyapunov function in $\TT_4$ by using \eref{e:prop}. To do so we specify the behavior of $\obar V$ on the boundary $\TT_4^*$.
Denoting throughout by $e_i$ the unit vector in direction $i\in \S$, by computing the two terms on the right hand side of \eref{e:prop} for $T_4$ as in \eref{e:t4} we obtain the following result.

\lmm{\label{l:v4} The function
\begin{equ}\label{e:v4}
  V_4(x) := {m_4^*} x_2^{\delta_4^*} + h_4 \sum_{k = 1}^{{x_1}} \frac{k^{\delta_4'}}{\delta_4''-2 + \binom{k}{5}5!} x_2^{\delta_4''-2}\,,
\end{equ}
approximating the solution of \eref{e:bvp} with
\begin{equ}\label{e:tt4}
h_4(x) := h_4 x_1^{\delta_4'}x_2^{\delta_4''} \qquad \text{and} \qquad V_4^*(x) := m_4^* x_2^{\delta_4^*}\,,
\end{equ}
for $h_4\in \RR$ and $m_4^* \in \RR$ is a well defined local Lyapunov function for all $x = (x_1, x_2) \in \TT_4$ and all $\delta_4^*, \delta_4'' > 0$, $\delta_4' \in \Rr$. Furthermore assuming that
\begin{equ} \label{e:hom4}
  \delta_4'' - 2 = \delta_4^*\,,
\end{equ}
we have $V_4 \ssim{w} l^{\delta_4}$ for $w = e_2$ and $\delta_4 := \delta_4^*$.}

\begin{proof}{Under the assumptions of this lemma we show in the appendix that the solution to \eref{e:bvp} is well-defined and can be approximated by \eref{e:v4}.
  We immediately notice that this function scales homogeneously under $\slw$ for $w = e_2$ iff $\delta_4''-2 = \delta_4^*$.
We now show that this function is a local Lyapunov function on $\TT_4$ for the generator $T_4$. Computing
  \begin{equs}
    T_4 V_4(x) &= x_2^2 \binom{x_1}{5}5!  \pc{V_4(x-5e_1) -V_4(x)} - x_2^3 \partial_2 V_4(x)\\
    & \leq - \binom{x_1}{5}5!  \frac{h_4}{\delta_4 + 5\binom{x_1}{5}5!} x_2^{\delta_4''} - \pc{m_4^* + \sum_{k = 1}^{{x_1}} \frac{h_4}{\delta_4 + \binom{k}{5}5!}}\delta_4 x_2^{\delta_4''} \leq - C (m_4^*+ h_4)\,,
  \end{equs}
  for a constant $C>0$.
}\end{proof}

\subsubsection{Transport regions: $\TT_1$, $\TT_2$ and $\TT_3$}
\label{s:transport}

Recall from \eref{e:transportgen} and \rref{r:asgenerator} that for $w \in S_+^1$ the scaled generator $\LL (\slw x)$ converges to a transport operators $\{T_i\}$. For this reason we call regions $\TT_1, \TT_2$ and $\TT_3$ \emph{transport regions}. In these regions, we obtain the solution to the Poisson equation \eref{e:bvp} by the method of characteristics \eref{e:prop}, as we do below.

We further recall that in $\LL$ is approximated by $T_i$ under the family of scalings $\slw$ for $w \in  \WW_i$. For this reason instead of defining \emph{one} scaling that maps the interior of the region $\TT_i$ to a compact we explore the asymptotic behavior of $\LL$ and the candidate local Lyapunov function $\obar V_i$ through a \emph{family} of scaling transformations. We use this fact to prove \eref{e:phi3} by the scaling analysis of Section~\ref{s:Lyapunov-2} carried out for all $w \in \WW_i$. To do so we need to assume that $h_i$ scales homogeneously under \emph{all} such transformations, \ie that there exists $\delta_i\colon S^1 \to \Rr$ such that
\begin{equ}\label{e:hscaling}
 h_i \ssim{w} l^{\delta_i'(w)} \qquad \forall\, w \in \WW_i\,.
\end{equ}
This condition is in particular satisfied by choosing $h_i$ to be a monomial, as we do below. We start by $\TT_3$.
\lmm{ \label{l:v3} The Lyapunov function $\obar V_3$ solving \eref{e:bvp} with
\begin{equ}\label{e:tt3}
h_3(x) := h_3 x_1^{\delta_3'}x_2^{\delta_3''}  \qquad \text{and} \qquad V_4(x) = m_4 x_2^{\delta_4}\,,
\end{equ}
for $h_3 > 0$ and $m_4 = m_4(b_2) > 0$ is well defined for all $x =  (x_1, x_2) \in \TT_3$ and all $ \delta_4> 0$, $\delta_3'' \in \Rr$ and $\delta_3' \neq 4$\,. Furthermore, for the choice of constants
\begin{equ}\label{e:hom3}
  \delta_3' > 4\,, \qquad \delta_3''-2 = \delta_4 \qquad \text{and} \qquad h_3 =(\delta_3'-4)m_4{b_2^{-(\delta_3'-4)}}\,,
\end{equ}
we can write $\obar V_3 (x) = h_3 (\delta_3' - 4)x_1^{\delta_3'-4} x_2^{\delta_3''-2}$.
}
\begin{proof}
  By the method of characteristics we obtain $\obar V_3$ satisfying \eref{e:bvp} by integrating $h_3$ along the solutions of the set of ordinary differential equations $\dot x = T_3 x$. Recalling by \eref{e:transportgen} that such solutions are moving from $x_1(0)$ to $b_2$ on lines with $x_2(t) = x_2(0)$, by our choice \eref{e:tt3} of boundary condition on $\TT_{34}$ we obtain
  \begin{equ}
   \obar V_3(x) = m_4 x_2^{\delta_4} + \int_{b_2}^{x_1} h_3 z^{\delta_3'}x_2^{\delta_3''}  \frac{1}{z^5x_2^2} \dd z =  m_4 x_2^{\delta_4} -\frac{ h_3 }{\delta_3'-4} b_2^{\delta_3'-4}x_2^{\delta_3''-2} +\frac{ h_3 }{\delta_3'-4}x_1^{\delta_3'-4}x_2^{\delta_3''-2}   ~,
  \end{equ}
  where we have assumed that $\delta_3' > 4$.
  This function is clearly well defined in $\RR^2$ for all choices of parameters. Now we see that in order for $\obar V_3$ to scale homogeneously under $\slw$ for all $w \in \WW_3$ we need \eref{e:hom3} as $h_3 > 0$. This directly implies that $\obar V_3$ has the desired form.
\end{proof}

\rmk{The requirement of having a Lyapunov function $\obar V$ that scales homogeneously under \emph{all} transformations $\slw$ for $w \in \WW_i$ can be relaxed to having $\obar V$ scale homogeneously under $\slw$ for $w = (1,1)/\sqrt2$. This method allows to construct a larger family of Lyapunov functions. However, in this case attention must be paid not to construct candidate Lyapunov functions that diverge at the boundary. We carry out such an alternative construction in the regions $\TT_1, \TT_2, \TT_3$ in the appendix.}

Proceeding to construct the candidate local Lyapunov function in $\TT_2$, we notice that this region is invariant under $\slw$ for $w = (1,1)/\sqrt 2$. Defining throughout for any $x \in \TT_i$, $\pi_{ij}x$ as the projection of $x$ onto $\TT_{ij}$ along the characteristics of $T_i$ in $\TT_i$ we construct $\obar V_2$ assuming that $h_2(\,\cdot \,)$ scales homogeneously under that transformation:
\lmm{ \label{l:v2} The Lyapunov function $\obar V_2$ solving \eref{e:bvp} with
\begin{equ}\label{e:tt2}
h_2(x) := h_2 ( x_1 + 5 x_2)^{\delta_2'} \qquad \text{and} \qquad V_3(\pi_{23} x) = m_3 (x_1 + 5 x_2)^{\delta_3}\,,
\end{equ}
for $h_2 > 0$ is well defined for all $x = (x_1, x_2) \in \TT_2$ and all $\delta_2' \in \Rr$, $ \delta_3> 0$. Furthermore, for the choice of constants
\begin{equ}\label{e:hom2}
  \delta_2' = \delta_3' + \delta_3''\,,
\end{equ}
we have $\obar V_2 \ssim{w} l^{\delta_2}$ for $w = (1,1)/\sqrt 2$ and $\delta_2 := \delta_2'-6$\,.
}
\begin{proof}
We again find the solution to \eref{e:bvp} in $\TT_2$ by the method of characteristics. Denoting by $\gamma_3(x,\pi_{23})$ the path along the characteristic of $T_2$ starting at $x$ and ending at $\pi_{23}x$ and noting that $h_2(\,\cdot\,)$ defined in \eref{e:tt2} is constant on such a path we have
\begin{equ}\label{e:integralv2}
 \obar V_2(x) = V_3(\pi_{23}x) + h_2(x_1+5x_2)^{\delta_2'} \int_{\gamma_3(x,\pi_{23}x)} \frac{1}{z_1^5z_2^2} \dd z\,.
\end{equ}
Consequently, using that $\pi_{23} x = (x_1+5x_2)(1+5 b_1)^{-1}(1,b_1)$ we write the explicit result of the integral as
\begin{equ}\label{e:explicitv2}
  \obar V_2(x) = m_3(x_1+5x_2)^{\delta_3} + \frac{h_2}{12}(x_1+5x_2)^{\delta_2'-6} \frac{5 + b_1^{-1} }{1 - b_1^{-2}} (P(b_1) - P(x_2/x_1))\,,
\end{equ}
for $P(x) := -(12/x) + 3000 x + 7500 x^2 + 12500 x^3 + 9375 x^4 + 300 \log x$, $\delta_3 := \delta_3'+\delta_3''-6$ and $m_3 = m_3(b_1) := h_3(b_1)^{\delta_3''}(1+b_1)^{\delta_3}$. Since the difference term is scale--invariant under $\slw$ for $w = (1,1)$ we obtain the homogeneous scaling behavior of \eref{e:explicitv2} iff $\delta_3 = \delta_2'-6$, leading to \eref{e:hom2}.
\end{proof}
\begin{remark}
  The homogeneous scaling behavior of $\obar V_2$ in \eref{e:explicitv2} could also be derived without explicit integration of \eref{e:integralv2}. Indeed, since both the integrand and the path length of $\gamma$ scale homogeneously in $l$ under $\mathscr S_l^{(1,1)/\sqrt 2}$, one must have $\obar V_2 \ssim{w} l^{\delta_2'-7+1}$ for $w = (1,1)/\sqrt 2$ in $\TT_2$ (up to a constant that increases in the parameters $b_0, b_1$).
\end{remark}
We conclude by considering the region $\TT_1$. Similarly to the case of $\TT_3$ we construct a local Lyapunov function choosing a function $h_1$ that scales homogeneously under the family of scaling transformations $\slw$ for $w \in \WW_1$.
\lmm{\label{l:v1}The Lyapunov function $\obar V_1$ solving \eref{e:bvp} with
\begin{equ}
h_1(x) := h_1 x_1^{\delta_1'}x_2^{\delta_1''}    \qquad \text{and} \qquad \obar V_2(x) = m_2 x_1^{\delta_2}\,,
\end{equ}
for $h_1 > 0$ is well defined for all $x = (x_1,x_2) \in \TT_2$ and all $ \delta_2> 0$, $\delta_1', \delta_1'' \in \Rr$ with $\delta_1'' \neq 1$\,. Furthermore, for the choice of constants
\begin{equ}\label{e:hom1}
  \delta_1'' < 1\,, \qquad \delta_2' = \delta_1' + \delta_1'' \qquad \text{and} \qquad h_1 = m_2\, (1-\delta_1'')b_1^{\delta_1''-1}\,,
\end{equ}
we can write $\obar V_1 (x) = h_1 (1-\delta_1'')^{-1} x_1^{\delta_1'-5} x_2^{\delta_1''-1}$.
}
\begin{proof}
  We obtain the Lyapunov function by integrating along the characteristic lines of the transport operator $T_1$. Noting that these lines satisfy $x_1(t) = x_1(0)$ for all $t > 0$ we write
  \begin{equ}
   \obar V_1(x) = \obar V_2(\pi_{12} (x)) + h_1 \int_{x_1/b_1}^{x_2} x_1^{\delta_1'} y^{\delta_1''} \frac{1}{x_1^5 y^2} \dd y = m_2 x_1^{\delta_2} - \frac{h_1 }{1-\delta_1''} x_1^{\delta_1'-5} { ( x_1/b_1)^{\delta_1''-1} } + \frac{h_1 }{1-\delta_1''} x_1^{\delta_1'-5} {x_2^{\delta_1''-1}}\,, \label{e:v1}
  \end{equ}
  for
  \begin{equ}\label{e:m2}
    m_2(b_1) = m_3 (1+5/b_1)^{\delta_3} + \frac{h_2}{12}\frac{5 + b_1^{-1} }{1 - b_1^{-2}}(P(b_1)-P(b_1^{-1}))(1 + 5/b_1)^{\delta_2} > 0\,.
  \end{equ}
  We immediately recognize that the right hand side of \eref{e:v1} scales homogeneously under $\slw$ for $w \in \WW_1$ if \eref{e:hom1} holds, resulting in the desired definition for $\obar V_1(\,\cdot\,)$\,.
\end{proof}

\subsubsection{Diffusive region: $\TT_0, \TT_0'$}
\label{s:diffusive}

Because noise dominates the behavior of the process for small values of $x_2$, throughout we refer to $\TT_0$ as the \emph{diffusive} region.

First, we continue to $\TT_0$ the Lyapunov function generated in the previous sections $\{x_2 = b_0\}$ by propagation through \eref{e:prop} in the transition region $\TT_0' := \{x \in \mathbb N_0^2~:~x_2 \in [2,b_2]\}$. Here, expectations are taken with respect to the process $X_t$ with asymptotic generator
\begin{equ}\label{e:t0}
  T_0' f(x) := k_3 x_1^5 \binom{x_2}{2} (f(x_1, x_2 + 1)-f(x))\,.
\end{equ}
This generator approximates the rate $\Lambda_r(x)$ as a power in the components that diverge under the scaling $\slw$ for the chosen $w$ while leaving the binomial formulation in the components that are not affected by such scaling. We refer to such a generator as a \emph{semi-continuous} approximating generator.

\lmm{\label{l:v0'}  The function $V_0'$ defined in \eref{e:prop} with $h_0'(x), V_1(x)$ given by
\begin{equ}
h_0'(x) := h_0' x_1^{\delta_0''} \qquad \text{and} \qquad V_1(x) = m_1 x_1^{\delta_1}\,,
\end{equ}
for $h_0' > 0$ and $m_1 = m_1(b_1) > 0$ is well defined for all $x = (x_1,x_2) \in \TT_{0}'$ and all all $\delta_0'' \in \Rr$, $ \delta_1> 0$. Furthermore, for the choice of constants
\begin{equ}\label{e:hom0'}
  \delta_1 = \delta_0''-5\,,
\end{equ}
we have that $\obar V_0' (x) \ssim{w} l^{\delta_1}$ for $w = e_1$.}
\begin{proof}{The above result is obtained by directly applying \eref{e:prop} for the dominant generator $T_0'$ from \eref{e:t0}:
$$\obar V_0'(x) = m_1 x_1^{\delta_0} + h_0' x_1^{\delta_0''-5} \sum_{k = x_2}^{b_1} \binom{k}{2}^{-1}\,,$$
which scales homogeneously under $\slw$ for $w = e_1$ iff \eref{e:hom0'} holds.}\end{proof}

We now consider the region $\TT_0$. By our fundamental assumption, we choose $h_0$ and $\obar V_0'$ to scale homogeneously under the transformation that leaves the set $\TT_0$ invariant ($\slw$ for $w = e_1$), \ie
\begin{equ}\label{e:delta0'}
  \begin{cases} h_0(x) &:= h_0 x_1^{\delta_0'}\\\obar V_0'(x) &:= m_0' x_1^{\delta_0}\end{cases} \qquad \text{for}\quad \delta_0' \in \Rr\,, \delta_0 \in \RR\,, h_0,m_0 \in \RR\,.
\end{equ}
We then obtain the following result for the propagated local Lyapunov function in $\TT_0$:
\lmm{\label{l:v0} Let $V_0$ be the function defined on $\TT_{0}$ through \eref{e:prop} with $h_0(x), \obar V_0'(x)$ specified in \eref{e:delta0'}. Then, for \abbr{crn0} $V_0$ is well defined for all choices of $\delta_0 \in \RR$ and $\delta_0' \in \Rr$, while the same is true for \abbr{crn1} only if $\delta_0 < 1$ and $\delta_0' < 0$.
If the above conditions are satisfied and we choose $\delta_0 > 0$, there exists a decreasing, positive function $m_0\colon (1,2) \to \RR$, such that the function
\begin{equ}
  V_0(x) = 
  \begin{cases}
    m_0(x_2) x_1^{\delta_0} & \text{for } x_2 \in (1,2)\\
    m_{0}(1) (x_1+1)^{\delta_0 } + h_0 x_1^{\delta_0'} & \text{for } x_2 = 0
\end{cases}
\label{e:v0}\, ,
\end{equ}
is a Lyapunov function in $\TT_{0}$ setting $\delta_0' = \delta_0$ for \abbr{crn0} and $\delta_0' = \delta_0' + 1$ for \abbr{crn1}.}

\proof{see appendix.}

\begin{remark}
  The Lyapunov function \eref{e:v0} does not scale homogeneously as  $x_1 \rightarrow \infty$ when $x_2 = 0$.
The leading order dynamics \eref{e:taylor}, which simply moves back and forth between $x_2=0$ and $x_2=1$ without changing $x_1$, does not faithfully capture the governing dynamics as $x_1\rightarrow \infty$. Rather than performing a systematic singular perturbation analysis, we choose to keep the entire generator in this region. Therefore we must choose a Lyapunov fiction which captures the interplay between terms in the generator and hence can not scale homogeneously.
\end{remark}

Combining the continuity conditions presented in Lemmas~\ref{l:v4}--\ref{l:v0'}, we obtain the following relationships between the exponents of our Lyapunov functions:
\begin{equ}
 \delta_1' = \delta_0'' =  5 + \delta_0\,, \quad \delta_1'' = \delta_2' - 5 - \delta_0\,,  \quad \delta_2' = \delta_3' + \delta_3''\,, \quad \delta_3'' = \delta_4'' = \delta_4^* + 2\,.
\end{equ}
Furthermore, the condition on compact sublevels sets of the Lyapunov function reads:
\begin{equ}
 \delta_0 > 0\,,\qquad \delta_1'' < 1\,, \qquad \delta_2'  > 6 \,, \qquad \delta_3' > 4~.
\end{equ}
In particular, for any $\epsilon \in (0, \delta_0/2)$ the following choice of constants works for \abbr{crn1}:
\begin{equs}
 \delta_0 \in (0,1)\,, \quad \delta_0' = \delta_0\,, \quad \delta_1' = \delta_0'' =  5 + \delta_0 \,, \quad &\delta_1'' = 1 - \epsilon \,,\quad \delta_2' = 6 + \delta_0- \epsilon \,,\\ \delta_3'=4+ \epsilon\,, \quad \delta_3'' = 2 + \delta_0 - 2 \epsilon& \,, \quad \delta_4 = \delta_0 - 2 \epsilon \,. \label{e:deltas}
\end{equs}
We introduce the $\epsilon>0$ in the definition of $\delta_1'$ and $\delta_1''$ in \eref{e:deltas} to enforce the conditions on $\delta_1''$ and $\delta_3'$ stated in \lref{l:v1} and \lref{l:v3} respectively. In doing so, we ensure that $V_1= C -\int_1^{x_2} s^{-1 \pm \epsilon} \dd s = -x_2^{\pm \epsilon} +C$ scales as a power. When $\epsilon=0$, $V_1= -\int_1^{x_2} s^{-1} \dd s = \log x_2$. This logarithmic scale would complicate the analysis.  Having chosen $\epsilon>0$, the Lyapunov functions $V_i$ scale at least as $l^{\delta_0-2 \epsilon}$ under the relevant transformations.

\subsection{Verification of the stability condition} \label{s:Lyapunov1}

 We now proceed to prove that the local Lyapunov functions defined in the previous section satisfy asymptotically the boundary value problem for the full generator \eref{e:generator} as stated in \lref{l:local}.

\lmm{\label{l:local}Let $i$ index a region of phase space. Then there exists a compact set $\KK$ and a constant $\Ch > 0$ such that the Lyapunov pair $(V_i,h_i)$ satisfies
\begin{equ}\label{e:local}
  \LL V_i(x) \leq - \Ch h_i(x) \qquad \text{for all} \quad x \in \TT_i \cap \KK^c\,.
\end{equ}
}

To prove \lref{l:local}, we apply the full generator \eref{e:generator} to the candidate local Lyapunov function $V_i$ obtained in the previous section and show that the corrections to the leading order term $h_i$ are negligible for large $l$. We proceed by considering each region of phase space separately.

\subsubsection{Diffusive region: $\TT_0$}

In the diffusive region $\TT_0$, the full generator of the Markov process was used to construct the Lyapunov function in \lref{l:v0}, so \eref{e:local} holds by construction. In particular, for $x_2 = 0$ we have
\begin{equ}
  \LL V_0(x) = \Lambda_1(x) \Delta_1 V_0(x) = m_0(1) (x_1+1)^{\delta_0} - (m_0(1) (x_1+1)^{\delta_0} + h_0 x_1^{\delta_0'}) = - h_0 x_1^{\delta_0'}\,,
\end{equ}
as expected. Similarly, for $x_2 = 1$ we obtain for \abbr{crn1} that for $x_1$ large enough
\begin{equs}
  \LL V_0(x) & = \sum_{i = 1}^2\Lambda_i(x) \Delta_i V_0(x)  = m_0(2) (x_1+1)^{\delta_0} - m_0(1) x_1^{\delta_0} + x_1 (m_0(1) (x_1+1)^{\delta_0} + h_0 x_1^{\delta_0'} - m_0(1) x_1^{\delta_0})\\\label{e:firstlevelcheck}
  & \leq (m_0(2) - m_0(1) ) x_1^{\delta_0} + m_0(1) \delta_0 x_1^{\delta_0} + h_0 x_1^{\delta_0'+1}\,.
\end{equs}
For small enough $h_0> 0$, because $\delta_0' + 1 = \delta_0 < 1$ the right hand side of the above expression is negative upon choosing $m_0(2) \in (0, m_0(1) (1- \delta_0) - h_0)$\,.



We now proceed to establish \eref{e:local} in the region $ \TT_0'$ as a special case of the following result. To do so, we define the \emph{semi--continuous} approximation to the reaction rates $\Lambda_r$ in direction $w \in \partial S_+^1$ as follows:
\begin{equ}
  \lambda_{r}^w(x) := \kappa_r \prod_{i \in \S\setminus\PP_w} x_i^{(c_\inn^r)_i} \prod_{i \in  \PP_w} \binom{x_i}{(c_\inn^r)_i} (c_\inn^r)_i!\,,
\end{equ}
where we have defined $\PP_w := \{i \in \S~:~w_i = 0\}$\,.

\lmm{\label{l:localjump} For any $i \in \S$, $b, b' \in \Nn_0$ with $b' > b \geq \max_{\R} (c_\inn^r)_i$ let  $\TT := \{x \in \Nn_0^2~:~x_i \in (b, b')\}$. Then if $\R_{w} = \{r^*\}$, \eref{e:local} holds on $\TT$ for $l$ large enough. }

\begin{proof} To study the behavior of $\LL V$ under $\slw$ for $w = e_j$ we start by writing
\begin{equ}\label{e:diff}
 \Delta_r V(\slw (x)) = \pc{V(\slw(x) + c^r) - V(\slw(x) + c_i^re_i)} + \pc{V(\slw(x) + c_i^re_i) -V(\slw(x))}~.
\end{equ}
and proceed to consider the two difference terms separately.
We start by the second term, corresponding to a jump in direction $i \in \PP_w$, and write
\begin{equ} \label{e:boundondiff}
 V(x + c_i^r e_i) -V(x) = - \Ex{x}{\sum_{k = 0}^{c_i^r-1} h(x + k e_i ) \Delta \tau_{x + k e_i} } = - \sum_{k = 0}^{c_i^r-1}\frac{h(x + k e_i )}{\lambda_{r^*}^w(x + k e_i)}  \leq  c_i^r \max_{|k| <|c_i^r|}\pc{{\frac{h(x + k e_i)}{\lambda_{r^*}^w(x + k e_i)}}}~,
\end{equ}
where $r^* \in \R_\WW$ and $\Delta \tau_{x}$ is the exponentially distributed jumping time of $X_t$ at $x\in \Nn_0^d$. Similarly, in the case of $r^*$ using that $h(x) > 0$ we have
\begin{equ} \label{e:boundondiff*}
 V(x + c_i^{r^*} e_i) -V(x) =  - \sum_{k = 0}^{c_i^{r^*}-1}\frac{h(x + k e_i )}{\lambda_{r^*}^w(x + k e_i)}  \leq -{{\frac{h(x )}{\lambda_{r^*}^w(x )}}}~.
\end{equ}
At the same time by combining our homogeneous scaling assumption with Taylor theorem for the second term we have
\begin{equ}\label{e:diffL}
 V(\slw(x) + c^r) - V(\slw(x) + c_i^re_i) \leq l^{\delta-1}(\partial_j \obar V(x + c_i^re_i) + l^{-1} R_{j,j}^r(x + c_i^r e_i,l)) ~,
\end{equ}
where $\obar V$ is $C^2$ in the $j$ direction and we define the remainder $R_{i,j}^r(y,l) := \sup_{a \in (-1,1)}|\partial_i\partial_j \obar V(y +  l^{-1} ae_i)|$.
Furthermore, we note that there exists a constant $\Clambda^r > 0$ such that for all $r \in \R$ and all $x \in \Nn_0^d$ such that $x - c_\inn^r \geq 0$ componentwise we have
\begin{equ}\label{e:lambdalbub}
\Clambda^r \lambda_r(x) \leq \Lambda_r(x) \leq \lambda_r(x)\,,
\end{equ}
where the right inequality holds by definition while the left one results from the increasing character of $x^{-a}\binom{x}{a}$ in $x>a$ for $a>0$.

 Combining \eref{e:boundondiff}--\eref{e:lambdalbub} and using the boundedness of $V$ in any compact set we obtain, for any $x \in \pi \TT := \{x \in \TT~:~x_j = 1\}$,
\begin{align}
 \LL V(\slw(x)) & = \sum_{r \in \R} \Lambda_r(\slw(x)) \pc{V(\slw(x) + c^r)- V(\slw(x))} \notag\\
 & \leq \sum_{r \in \R} \Lambda_r(\slw(x)) \Big(V(x + c_i^r e_i) -V(x)
   + l^{\delta-1}(\partial_j \obar V(x + c_i^r e_i) + l^{-1} R_{j,j}^r(x + c_i^r e_i,l)) \Big)\notag \\
   &\leq - l^{\delta'(w)} \Bigg(  h(x) \pq{  C_r^\lambda - l^{-1}c^* \sum_{r \in \R\setminus \R_w} \frac{\Lambda_r(x)}{h(x)}\max_{|k| <|c_i^r|}\frac{h(x + k e_i)}{\lambda_{r^*}^w(x+ke_i)}}
  + l^{-1}\pc{ \sum_{r \in \R} \Lambda_r(x)\Cv } \Bigg) \,,\raisetag{-.5em}\label{e:final}
\end{align}
where in the second inequality we have used that $h(x) \ssim{w} l^{\delta+\dtp{c_\inn^{r^*},w} }$, that $\Lambda_r(\slw x)/\lambda_{r^*}^w(\slw x) \leq l^{-1}\Lambda_r( x)/\lambda_{r^*}^w( x)$ for all $r \in \R\setminus \R_w$ if $w \in \{e_i\}$ and we have bounded from above the derivative terms
$\partial_j \obar V(x + c_i^re_i) + l^{-1} R_{j,j}^r(x + c_i^r e_i) \leq \Cv$ on the compact $\pi \TT$.
The boundedness of the $x$-dependent term in square brackets on the right hand side of \eref{e:final} on the finite set $\pi \TT$ and the divergence of $l$ prove the desired result.
\end{proof}

\subsubsection{Transport regions: $\TT_1$, $\TT_2$ and $\TT_3$}

As anticipated in the previous section, the function $\obar V_i$ on $\TT_1$ (resp. $\TT_3$) is assumed to scale homogeneously under $\slw$ for all $w \in \WW_1$ (resp. $\WW_3$). We use this fact to explore the asymptotic behavior of $\obar V_i$ by writing every point $z \in \RR^d$ in terms of \emph{toric coordinates}   $\pc{\theta,w} : \pc{\RR^{n}}^o \rightarrow  \Rr_{>1} \times S^{n-1}$ defined by
\begin{align}
 \theta(z) := \exp(\|\log((2c^*)^{-1}z)\|_2)\,, \quad
w(z) :&=\frac{1}{\log \theta(z)} \log((2c^*)^{-1}z) \quad \text{where} \quad z = [\theta(z)^{w(z)_i}]_{i \in \S} \label{e:tropical}
\end{align}
with $c^* := \sup_{r \in \R}\{\|c_\inn^r\|_1,\|c_\outt^r\|_1\}$ and $\log\colon\RR^d \to \Rr^d$ represents the component-wise logarithm. This transformation maps a $c^*$ neighborhood of any point $x \in \Nn_0^d$ such that $x_i > 2c^*$ to the $c^*$-neighborhood of the point $z^*$ with $z_i^* := c^*$ for all $i \in \S$.

Under the assumption of homogeneous scaling of the Lyapunov function, the generator of the Markov process can be asymptotically approximated along toric rays in $S_+^1$ by the generator of the transport process from \eref{e:Ltransport}. This convergence happens pointwise in $w$, and is therefore not sufficient for ensuring the required scaling property, which must hold \emph{uniformly} in $w$. This uniform convergence is established in \lref{l:taylor2} below.

\lmm{\label{l:taylor2} For all $\epsilon > 0$ and $r \in \R$ there exist $\Cr \in (1-\epsilon, 1 + \epsilon)$ and $\Cx > 0$ such that for all $w \in S_+$ we have, for $l$ large enough,
 \begin{equs}
  \LL V(\slw(x)) \le   \sum_{r \in \R}l^{\dtp{c_\inn^r, w}+v_r(w)}  \Cr  \lambda_r(x)  \dtp{ (\slw)^{-1}(c^r),\nabla \obar V(x)}\qquad \text{for all } x \in \Nn_0^d,\, x_i > \Cx \, \forall i \in \S\,.
  \label{e:taylor2}
 \end{equs}
}
\begin{proof}
 Consider the scaling of the generator \eref{e:generator} along an arbitrary toric ray. By definition of the transformation in \eref{e:tropical} it is sufficient to know the value of $\obar V$ on a compact set $\KK^*$ to obtain through a scaling transformation the values of $\obar V$ outside of $\KK^*$.

 Using the homogeneity of $V$ we bound the differences $\Delta_r$ from \eref{e:diff} similarly to \eref{e:diffL}:
 \begin{equs}
   V(\slw(x) + c_i^r e_i) -V(\slw(x)) & \leq l^{v_r-w_i}(\partial_i \obar V(x ) + l^{-w_i} R_{i,i}^r(x ))\,,\\
  V(\slw(x) + c^r) - V(\slw(x) + c_i^re_i) &\leq l^{v_r-w_j}(\partial_j \obar V(x + l^{-w_i}c_i^re_i) + l^{-w_j} R_{j,j}^r(x + l^{-w_i}c_i^r e_i))
  \\&\leq l^{v_r-w_j}(\partial_j \obar V(x) + l^{-w_i}R_{i,j}^r(x) +  l^{-w_j} R_{j,j}^r(x + l^{-w_i}c_i^r e_i)) ~,
\end{equs}
where we have chosen the indices $i, j \in \S$ such that $w_j\geq w_i>0$.
  Note that by the boundedness of $\obar V$ and its partial derivatives, we have that $\sup_{\KK^*} (R_{i,i}^r + R_{j,j}^r + R_{i,j}^r) \leq K$ for all $r \in \R$.
  Consequently, choosing $\Cx$ large enough for $l^{-w_i} < \epsilon/6K$ for all $i \in \S$ we have
 \begin{equ}\label{e:diffL3}
   \Delta_r (V \circ\slw(x)) \leq \Clambda' \pc{l^{-w_1}c_1^r \partial_1 \obar V(x) + l^{-w_2} c_2^r\partial_2 \obar V(x)}\,,
 \end{equ}
 where $\Clambda' \in (1-\epsilon/3, 1 + \epsilon/3)$. Because of $\lim_{x \to \infty}x^{-a}\binom{x}{a}a!=1$ we can choose the constant $\Clambda^r \in (1 - \epsilon/3, 1)$ in \eref{e:lambdalbub} upon possibly increasing $\Cx$ further, and we finally obtain the desired result by combining \eref{e:lambdalbub} with \eref{e:diffL3} and upon choosing $\Cr := \Clambda' \cdot \Clambda^r$.
\end{proof}

We now use the above result to prove that the candidate Lyapunov function satisfies \eref{e:local}. Denoting by $\TT_{ij}(a) := \{x \in \RR^2~:~\inf_{y \in \TT_{ij}}\|x-y\|_2 < a\}$ we use approximation \eref{e:taylor2} for $\slw(x) \in \TT_i \setminus \bigcup_j \TT_{ij}(\Cx)$ for $\Cx$ large enough and we obtain
\begin{equs}\label{e:firstbound}
 \LL \obar V_i(\slw(x)) & \leq \sum_{r \in \R } l^{v_r(w) } \sum_{i \in \S} l^{\dtp{c_\inn^r - e_i ,w}} \Cr \lambda_{r}(x) c_{i}^{r}\partial_i \obar V(x) ~.
\end{equs}
We note that the dynamics associated to the operator on the right hand side are, up to a change of constants $\kappa_r$, the mass action ordinary differential equations for the \abbr{crn} $(\S, \R', \C')$ with
\begin{equ}\label{e:WW}\R' := \{(r,i)\in \R \times \S ~:~ c^r_i \neq 0 \} \quad\text{and} \quad \C_\inn':=\{c_\inn^{r,i}~:~(r,i) \in \R'\}\,, \end{equ}
where we define $c_\inn^{r,i} := c_\inn^r-e_i$ and to each $(r,i) \in \R'$ we associate a reaction vector $c^{r,i} := c_i^r e_i$. Furthermore, we define for $w \in S_+^1$ the set of reactions that are \emph{exposed} by all $w \in \WW$ as
\begin{equ} \label{e:exposed}
  \R_{\WW}' := \{(r,i) \in \R'~:~\dtp{c_\inn^{r,i},w}\geq\dtp{c_\inn^{r',i'},w}\,\forall\, (r',i') \in \R', w \in \WW \}\,.
\end{equ}
We recognize that the regions $\WW_{0, \dots, 4}$ correspond to the partition $\WW^*$ associated to the convex hull $\WW$ of points in $\C_\inn'$ as defined in \cite{agazzi18} and proceed to apply \cite[Lemma~4.26 (d)]{agazzi18} to the present framework. To map this problem to the one in \cite{agazzi18} we make the change of notation $\R(\PP)_+ \to \R_{\WW_j'}$ where $w \in \WW$, $\R(\PP)_- \to \R' \setminus \R_\WW'$ and $\theta \to l$, so that $K \epsilon(\theta) \to \Cx'/\log l$. Doing so we have that for all $w \in \WW(\Cx')$
\footnote{recall that the set \eref{e:exposed} is nonempty for an $\WW$ by \cite[Remark~4.21]{agazzi18}.}
 \begin{equ}\label{e:l426}
   \dtp{c_\inn^{r',i'}- c_\inn^{r,i},w} \leq - \Cx'/\log l \qquad \text{for all } (r,i) \in \R_{\WW_j}', (r',i') \in \R'\setminus \R_{\WW_j}'\,.
 \end{equ}
Using \eref{e:l426} we bound the exponent of the scaling parameter for the subdominant terms in \eref{e:firstbound} obtaining that for any $\epsilon > 0$ we can choose $\Cx'$ large enough such that
\begin{equ}
  \lambda_{r'}(\slw(x)) \leq \epsilon/m \lambda_r(\slw(x)) \qquad \text{for all } \quad r \in \R_{\WW_j}\,,r' \in \R\setminus \R_{\WW_j}\,.
\end{equ}
 Now, using that $v_r(w) = \delta(w)-\min_{\S}w_i$ for all $r \in \R$ and that the derivatives of $\obar V$ are bounded away from $0$ in $\KK^*$, we see that for all $w \in \WW_j \setminus \bigcup_k \WW_{k}(\Cx')$ there exist for all $r \in \R$ constants $\Cr'(\epsilon) \in (1-2\epsilon,1+2 \epsilon)$ such that we have
\begin{equ}
  \LL V(\slw(x))
  \leq \sum_{(r,i) \in \R_{\WW_j}'} l^{\dtp{c_\inn^{r,i} ,w}+ v_r(w) } \Cr'(\epsilon) \lambda_{r}(x) c_{i}^{r} \partial_i \obar V(x)\,.
\end{equ}
Recalling definition \eref{e:deltass}, by continuity in $\epsilon$ of the right hand side and knowing that \eref{e:local} holds for $\Ch = 1$ in the limit $\epsilon \to 0$ we obtain desired result \eref{e:local} for any $\Ch \in (0,1)$ upon choosing $\epsilon$ small enough.

\subsubsection{Priming region: $\TT_4$}

Recalling that the propagated Lyapunov function in this region approximates the solution to \eref{e:bvp} with the ansatz \eref{e:hom4} for the leading order generator $T_4 f = \binom{x_1}{5}x_2^2 (f(x - 5 e_1) - f(x)) +x_2^3 (f(x - e_2) - f(x))$, we have for $w = e_2$ that
\begin{equ}\label{e:v4local}
  \LL V_4( \slw (x)) =  T_4 V_4( \slw(x))) + (\LL - T_4) V_4( \slw(x)) = l^{\delta_4''} \pq{T_4 \obar V_4(x) + l^{-\delta_4''} (\LL - T_4) V_4( \slw(x))}\,.
\end{equ}
We proceed to show that the second term in the square brackets goes to $0$ as $l \to \infty$. By \eref{e:boundondiff} and \eref{e:diffL} we can approximate (up to a multiplicative constant $\Cr > 0$) for large $x_2$ the difference terms in direction $x_2$ with partial derivatives, and obtain
\begin{equs}
(\LL-T_4) V_4(x) & = \sum_{r_1,r_2} \Lambda_r(x) (V_4(x + c^r)- V_4(x))\notag\\&\leq - \Cr x_2 \partial_2 \obar V_4(x) +  \Cr \partial_2 \obar V_4(x + e_1) + (\obar V_4(x + e_1) -\obar V_4(x)) + \indicator_{\{x_1\geq 5\}} \Cr x_2^2 \partial_2 \obar V_4(x - 5e_1 )\notag
\\& \leq  x_2^{\delta_4+1} \Cr {\delta_4} h_4 \pc{\indicator_{\{x_1\geq 5\}}B_5(x_1) G(x_1-5) + \frac{G(x_1+1)}{ x_2^2} + \frac{G(x_1)}{x_2}+  \frac{1 }{x_2\delta_4(\delta_4 + B_5(x_1+1))}  }~.
\end{equs}
where $G(x_1) := \sum_{i = 1}^{x_1} (B_5(i)+(\delta_4''-2) )^{-1}$ and $B_5(i) := \binom{i}{5} 5!$ and we used that $\delta_4 + 2 = \delta_4''$ from \eref{e:hom4}. Finally, we bound the right hand side of the above expression by $C B_5(x_1) G(x_1) x_2^{\delta_4+1}$ for a large enough $C>0$ and we obtain that the second term in \eref{e:v4local} scales as $l^{-1}$, proving the claim.


\section{Assembling a global Lyapunov function} \label{s:Lyapunov2}

This section is devoted to verifying that the local Lyapunov functions defined in the previous sections can be assembled
to generate a \emph{global} Lyapunov function. To guarantee that this is the case, we show that at the interface between two contiguous nonoverlapping regions of phase space the application of the generator to the global candidate global Lyapunov function is indeed negative and scales as required.
 We prove this in three steps. First, we introduce an intuitive condition for the assembly between two regions with candidate local Lyapunov functions solving \eref{e:bvp} to be a global Lyapunov function on the union of such domains. Then, we show that for any $h_j$ there exists a choice of parameters $h_i > 0$ such that the assembly works automatically for all the interfaces. Finally, we show that the choice of parameters above does not affect the relevant properties of the candidate global Lyapunov function far from the patching boundary.

\subsection{A condition of natural assembly}

In a recent series of papers \cite{mattingly151,mattingly152}, the problem of assembling local, homogeneously scaling Lyapunov functions was studied to prove stability of a certain family of diffusion processes. In that paper the authors show, under geometric assumptions related to the convexity of the (continuously assembled) Lyapunov function across the boundary separating two contiguous regions of phase space, that the assembled Lyapunov function automatically satisfies the desired Foster-Lyapunov condition on the union of those regions, and in particular on their common boundary. We will refer to these conditions as \emph{the interface curvature} condition. It ensures that the term, analogous to the \emph{Tanaka} or flux term derived in \cite{peskir07}, have the properties needed to avoid the often less intuitive smoothing/mollification procedures used in this assembly process. In the continuous diffusions setting, this flux term arises in the generalized It\^o formula, called \emph{Tanaka's formula}, because the Lyapunov function is only $C^1$ along the interface rather than the usual $C^2$ required by It\^o's formula.

\begin{figure}
 \centering
 \def\svgwidth{.30\textwidth}
 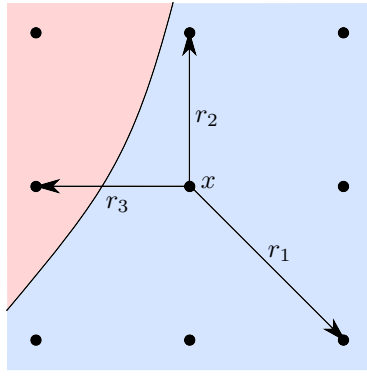
 \caption{Representation of the possible jumps of the process at $x \in \Rr^2$ close to the boundary between two regions with different local Lyapunov functions (in red and blue). In the generator $\LL$, only summands corresponding to jumps across the boundary ($r_3$) will feel the difference, while all the other terms will not ($r_1$ and $r_2$).}
 \label{f:tanaka}
\end{figure}

In this section, we adapt such conditions for ``natural assembly'' of local Lyapunov functions to the discrete setting.
To establish the interface curvature condition, we study the behavior of $\LL V$ close to an interface between two neighboring regions $\TT_i$ and $\TT_j$ with respective local Lyapunov function $V_i$ and $V_j$ to identify the equivalent of Tanaka's term in our setting. We take $V$ to be the (asymptotically) continuous assembly of $V_i$ and $V_j$ along $\TT_{ij}$. We
immediately note that, in general, only those terms corresponding to jumps across $\TT_{ij}$ will feel the discontinuity, as represented in \fref{f:tanaka}. We refer to such terms as cross-terms, and denote the corresponding reactions as $\R_{c}(x)$.  
Then, assuming without loss of generality that $x \in \TT_i$ we can write
\begin{equs}
 \LL V(x) & = \sum_{r \in \R\setminus \R_c(x)} \Lambda_r(x) \pc{V_i(x + c^r)-V_i(x)} + \sum_{r \in \R_c(x)}\Lambda_r(x) \pc{V_j(x + c^r)-V_i(x)}\\&
 = \LL V_i(x) + \sum_{r \in \R_c(x)}\Lambda_r(x) \pq{V_j(x + c^r) - V_i(x+c^r)} \\&
 = \LL V_i(x) + \sum_{r \in \R_c(x)}\Lambda_r(x) \pq{(\obar V_j(x + c^r) - \obar V_j(x + \beta_r(x)c^r)) - (\obar V_i(x + c^r)-\obar V_i(x+\beta_r(x)c^r))}~,\label{e:tanaka}
\end{equs}
where in the last equation, for every $r \in \R_c(x)$ we have chosen $\beta_r(x) \in (0,1)$ such that $x+\beta_r(x) c^r \in \TT_{ij}$ and we have used the continuity of $\obar V$. Throughout, we will refer to term in square brackets in \eref{e:tanaka} as the \emph{Tanaka or Flux term}
$\FF_V^r(x)$ and define the operator $\nabla_{x,r} V_i := (\obar V_i(x + c^r) - \obar V_i(x + \beta_r(x)c^r))$~.

To prove  the result of this section, we now note that $\LL V(x)$ satisfies \eref{e:bvp} along $\TT_{ij}$ if
\begin{axx}[Interface Flux Condition]  \label{a:tanaka}
 For all $x \in \mathbb N_0^d$ such that $\R_c(x) \neq \emptyset$, for any $r \in \R_c(x)$ one of the two following conditions holds:
 \begin{myenum}
  \item the Tanaka term is of negative sign \ie $\FF_V^r(x)<0$, or
  \item it is dominated in absolute value by $\LL V$ in the scaling that leaves $\TT_{ij}$ invariant.
 \end{myenum}
\end{axx}
\noindent Throughout we refer to the condition above as \emph{interface flux condition}. It ensures that the term generated by any jump across the boundary either has the correct sign to be comparable with direction of the desired inequality if neglected or that it scales in such a way to be negligible when compared to the dominant term along the relevant paths to infinity.

Condition~\ref{a:tanaka} (a) can be directly verified in a simple and geometrically intuitive way considering the convexity property of the Lyapunov function across the interface. To introduce the notion of $\bar V$ on a boundary $\TT_{ij}$ we define throughout the function $\TT_{ij}~:~\RR^d \to \Rr$ whose zero-level set is the surface $\TT_{ij}$ \ie $\TT_{ij} = \{x \in \RR^d~:~\TT_{ij}(x) = 0\}$\,.
\dfn{[Discrete Interface Curvature] \label{r:tanaka} We define the scalar curvature of the Lyapunov function $\obar V$ across the boundary $\TT_{ij}$ as
 \begin{equ}\label{e:curvature}
   \kappa(\obar V,x) := \lim_{\epsilon \to 0}\langle{c_\perp(x), { \nabla \obar V(x+ \epsilon c_\perp )- \nabla \obar V(x- \epsilon c_\perp )}\rangle}\,,
 \end{equ}
 where $c_\perp(x) :=  \nabla \TT_{ij}(x)/\|\nabla \TT_{ij}(x)\|_2$ and $\TT_{ij}$ is invariant under $\slw$ for a $w \in S_+^1$.
 For interfaces $\TT_{ij}$ that are parallel to $\partial \RR^2$, \ie that are invariant under $\slw$ for $w \in \{e_1,e_2\}$ we define for $|\alpha| \leq c^*$\, the discrete curvature of $\obar V$ at $x \in \TT_{ij}$ as
  \begin{equ}
    \label{e:curvature2}
    \obar \kappa_\alpha(\obar V,x) := {{  \obar V_i(x+ \alpha c_\perp(x) )-\obar V_j(x+ \alpha c_\perp(x) )}}\,.
  \end{equ}}
Now, if the curvature \eref{e:curvature} resp.~\eref{e:curvature2} is negative, the value of the function $\obar V_j$ in a point across but close enough to the interface is smaller than the one of the Lyapunov function $\obar V_i$ continued analytically to the same point. As shown in \lref{l:curvature} below, this makes the term $\FF_V^r(x)$ negative and automatically verifies Condition~\ref{a:tanaka} (a).

\lmm{[Discrete Interface Curvature Condition] Let $\TT_{ij}$ be invariant under $\slw$ for $w \in \obar S_+^1$. If $\kappa(V,x) < 0$ for all $x \in \TT_{ij}$ [respectively $\obar \kappa_\alpha(V,x) < 0$ for all $x \in \TT_{ij}(c^*)$, $|\alpha| < c^*$] with $\|x\|_2$ large enough, then $\FF_V^r(x)< 0$ for all $x \in \TT_{ij}$\,.
\label{l:curvature}}
\begin{proof}
  We start by proving the desired result for $\kappa(\obar V,x)$, \ie if $w \in S_+^1$. By \lref{l:taylor2} we can write for all $r \in \R_c(x)$ that $\nabla_{r,x}\obar V_j (x) = \dtp{c^r, \nabla \obar V_j(x_r^*)} + \epsilon$ for $\epsilon$ small at will provided that $x_r^* := x + \beta_r(x)c^r \in \TT_{ij}$ is large enough. Then, defining $c_\perp^r(x) := c_\perp(x_r^*) \dtp{c_\perp(x_r^*),c^r}$ we write
\begin{equs}
  \FF_V^r(x) = {\dtp{c^r, \nabla \obar V_j(x_r^*) - \nabla \obar V_i(x_r^*)}} =  {\dtp{c^r_{\perp}(x_r^*),\nabla \obar V_j(x_r^*)-\nabla \obar V_i(x_r^*)} + \dtp{c^r - c^r_{\perp}(x_r^*),\nabla \obar V_j(x_r^*)-\nabla \obar V_i(x_r^*)}} \,.
\end{equs}
By continuous differentiability of $\obar V$ along the boundary we note that the second term on the right hand side vanishes and we obtain the desired result by identifying the first summand with \eref{e:curvature}.

We now proceed to consider $\obar \kappa_\alpha(\obar V,x)$. Noting that in this case $x + c^r - c_\perp^r(x_r^*) \in \TT_{ij}$ we obtain
\begin{equs}
  \FF_V^r(x) &  = \obar V_i (x + c^r) - \obar V_j (x + c^r)  =  \obar V_i((x + c^r - c_\perp^r(x_r^*)) + c_\perp^r(x_r^*)) - \obar V_j((x + c^r - c_\perp^r(x_r^*)) + c_\perp^r(x_r^*))\,.
\end{equs}
Identifying the right hand side of the above equation with $ \obar \kappa_{\alpha}(\obar V,x + c^r - c_\perp^r(x_r^*))$ for $\alpha = |c_\perp(x^*)|/(1-\beta_r(x))$ completes the proof.
\end{proof}

\rmk{We recall that the negativity of the curvature \eref{e:curvature} and \eref{e:curvature2} is related to the existence of super- and sub-solution to the partial differential equation \eref{e:bvp} across the boundary $\TT_{ij}$.}

 \subsection{Scaling at the boundary}

Assuming that the constructed Lyapunov function is continuous at the boundary, use \rref{r:tanaka} to avoid the customarily lengthy calculations needed to assemble local Lyapunov functions between different definition domains. The idea relies on the tuning of the parameter $h_i  > 0$ in the regions at the interface in order to make the curvature conditions \eref{e:curvature}, \eref{e:curvature2} in  \rref{r:tanaka} automatically verified. In our case, however, this cannot be done in full generality, as the following example shows.

\exm{\label{ex:counter} We consider the partition of phase space for the networks \eref{e:e0}--\eref{e:e2} with boundaries \eref{e:boundaries} and study the convexity of the assembly of $ V_2 $ and $ V_1$ in a neighborhood of $\TT_{12}$. To do so, by our homogeneous scaling assumption and by the assumed continuity of the assembled function $V$ along $\TT_{12}$, it is sufficient to consider \eref{e:curvature} for $c_{\perp} = - c^{r_3}$.

In this case, for $x \in \TT_{12}$ we have by \eref{e:prop} and \eref{e:hom3} that
\begin{equ}
   \dtp{c^{r_3},  \nabla\obar V_{2}(x)} = \dtp{c^{r_3}, \nabla_y {\int_{\pi_{12} (x+y)}^{x+y} \frac{h_2(z)}{{\lambda_3(z)}} \dd z}} = -\|c^{r_3}\|_1 \frac{h_2(x)}{\lambda_3(x)} = -\|c^{r_3}\|_1 \frac{h_2 (x_1 + 5x_2 )^{\delta_2'}}{x_1^5x_2^2}\,,\label{e:counter1}
\end{equ}
where we recall that $\pi_{ij}$ is the projection of points in $\TT_i$ onto $\TT_{ij}$ along the characteristic lines used to construct $h_i$ and in the last equality we have used that the characteristics in $\TT_2$ are parallel to $c^{r_3}$. Similarly, for $V_1$ we have
\begin{equ}
  \dtp{c^{r_3} , \nabla \obar V_{1}(x)} = \dtp{c^{r_3}, \nabla \pq{h_1 x_1^{\delta_1'-5}x_2^{\delta_1''-1} }} = - h_1\pc{5\frac{\delta_1'-5}{x_1} + \frac{1-\delta_1''}{x_2}} x_1^{\delta_1' - 5} x_2^{\delta_1'' - 1} \,.\label{e:counter2}
\end{equ}
In particular for $x \in \TT_{12}$, \ie for $x_1 = b_1 x_2$ we have for \eref{e:counter1} and \eref{e:counter2}, respectively
\begin{equ}\label{e:crnablat12}
  \dtp{c^{r_3}, \nabla\obar V_{2}(x)} = -\|c^{r_3}\|_1 h_2 (5/b_1+1)^{\delta_2'}b_1^{2} x_1^{\delta_2'-7}\,,
  \quad
  \dtp{c^{r_3}, \nabla \obar V_{1}(x)} = - h_1\pc{5({\delta_1'-5}) + ({1-\delta_1''}){b_1}} b_1^{1-\delta_1''} x_1^{\delta_1' + \delta_1'' - 7} \,.
\end{equ}
Now, recalling from \eref{e:hom1} that $h_1 = m_2 (1-\delta_1'') b_1^{\delta_1''-1} $, we combine this with the expression \eref{e:m2} for $m_2$ to obtain the upper bound
$$\dtp{c^{r_3}, \nabla \obar V_{1}(x)} \leq - m_2'(b_1) h_2 x_1^{\delta_2' - 7} \qquad \text{for }\quad  m_2'(b_1) :=   \frac{({1-\delta_1''})^2(b_1^{-1} + 5 ) }{12  (1 - b_1^{-2}) }(P(b_1)-P(b_1^{-1})) b_1\,.$$
Combining the estimates from above we have that
\begin{equ}
  x_1^{7-\delta_2'}\dtp{c^{r_3} , \nabla\obar V_{2}(x) -  \nabla \obar V_{1}(x)} \geq - \|c^{r_3}\|_1h_2 (5/b_1+1)^{\delta_2'} b_1^{2} +    h_2 m_2'(b_1)  = h_2\pc{     m_2'(b_1)-\|c^{r_3}\|_1 (5/b_1+1)^{\delta_2'} b_1^{2}} \,,
\end{equ}
We see that the sign of above expression, which bounds the sign of the curvature $\kappa(\obar V,x)$ from below, is positive for large $b_1$ and independent of the parameter $h_2>0$. This parameter can therefore not be used to correct the curvature of $V$  at this interface. Because the parameter $b_1 $ was already bounded from below in the previous section, and that $\lim_{b_1 \to \infty}m_2'(b_1)/b_1^2 = \infty$ we see that there is no choice of parameters $b_i$, $h_i$ s.t. \eref{e:curvature} is satisfied in the general case.
\label{ex:1}}

We circumvent the problem highlighted in \exref{ex:1} by introducing a construction to tune the curvature of $\obar V$ on the boundary of interest while affecting only marginally its value on that set. This is possible because the value of $\obar V$ at the interface is obtained by integration of $h_i(x)$ in \eref{e:prop} across the set $\TT_i$ while the gradient $\nabla \obar V$ is a strictly local quantity, \ie it only depends on the choice of $h_i(\,\cdot\,)$ close to the boundary. In particular, changing the value of $h_i(x)$ in a small enough neighborhood of $\TT_{ij}$ will change the value of $\nabla \obar V$ but will have little effect on $\obar V$. We perform this change along the characteristics of the asymptotic generator $T_i$ and in a ``smooth'' way, \ie over many small enough steps, in order for \eref{e:local} to hold at the interface between regions where $h_i(x)$ is varied.

To realize the program outlined above (in this particular case of $\TT_2$), we dissect the problematic region $\TT_{i}$ into $n_i$ nonoverlapping radial subsets $\{\TT_i^{(k)}\}_{k \in (1, \dots, n_i)}$ with $\obar \TT_i^{(j)} \cap \obar \TT_i^{(j')} = \emptyset$ if $|j-j'| > 1$. In each cone, we then define the local Lyapunov functions $\obar V_i^{(k)}$ by \eref{e:prop} with
\begin{equ}
  h_i^{(j)}(x) := (\eta_i^*)^{j} h_i(x)\,,
\end{equ}
  for $\eta_i^* > 1$ and all $j>0$.  The procedure outlined above constructs a discrete interpolation given by the set-function pairs $\{\TT_i^{(j)},\obar V_i^{(j)}\}_{j \in 1 , \dots , n_i}$. We denote the new candidate local Lyapunov function assembled over the union of all $\TT_i^{(j)}$ by $\tilde V_i(\,\cdot\,)$, namely
  \begin{equ}
    \tilde V_i(x) := V_i^{(k)}(x) \qquad \text{for } x \in \obar \TT_{i}^{(k)}\,.
  \end{equ}

\begin{figure}
 \centering
 \def\svgwidth{.54\textwidth}
 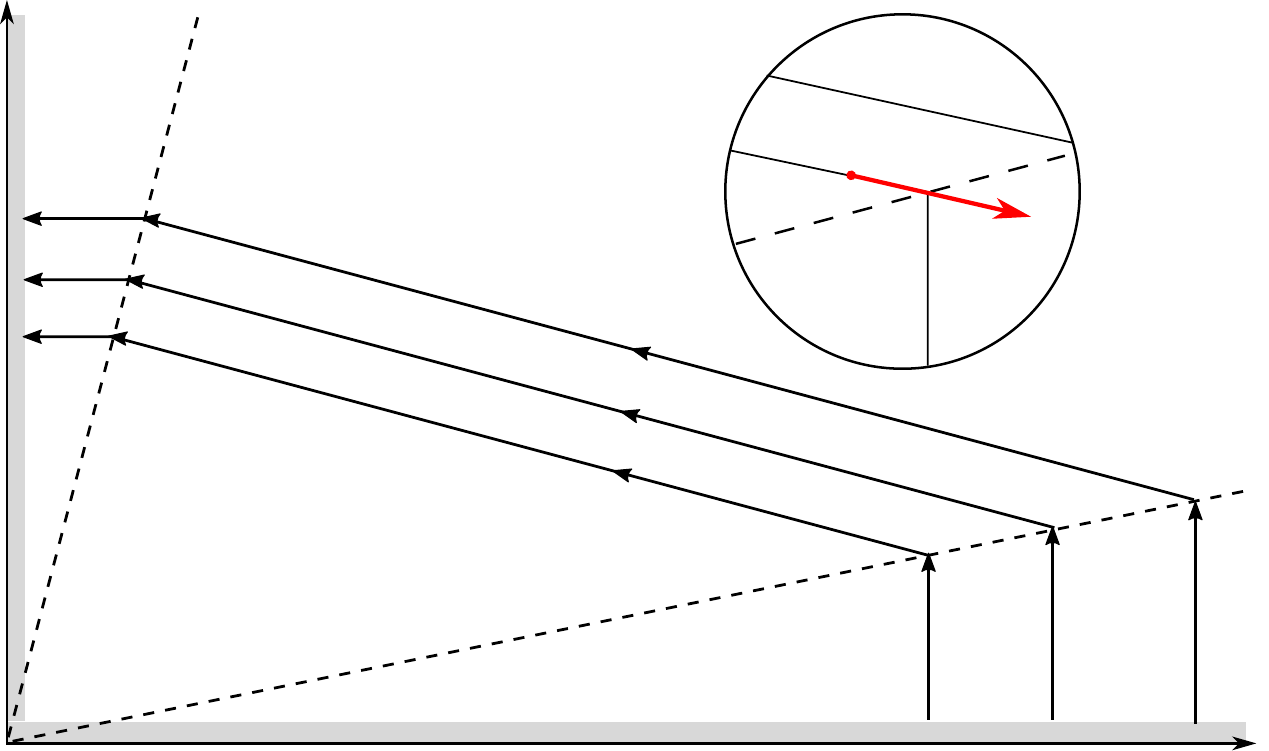
 \caption{Representation of the characteristic lines (solid vectors) of the dominant generators in the transport regions. The regions with different dominant generators are separated by dashed lines: Black dashed lines represent the original separation of dominance regions $\TT_{ij}$ while the green ones arise from the construction of $\tilde V_i$. The grey regions represent the diffusive and the priming regions. In the inset, the direction of the derivative in \eref{e:counter1} across the boundary $\TT_{12}$ is represented as a red vector. The parameter $\alpha_{}$ defines the slope of $\TT_{12}$ and $\TT_{23}$, so $\cos \alpha = b_1^{-1}$.}
 \label{f:longpatching}
\end{figure}

\lmm{\label{l:tuning} Assume that $\TT_{i}$ is invariant under $\slw$ for a $w \in \WW_i$ and that $\R_{\WW_i} = \{r^*\}$. Then for any $\epsilon > 0$, $M > 1$ in the partition $\{\TT_i^{(j)},\obar V_i^{(j)}\}_{j \in 1 , \dots , n_i}$ constructed above, the value of $n_i$, $\eta_i^*$ and the sets $\TT_{i}^{(j)}$  can be chosen so that $\tilde V_i$ is a Lyapunov function on $\TT_i$ and  for all  $x \in \TT_{ij}$ large enough, one has
\begin{equ}\label{e:disentangle}
  |\obar V_i(x)-\tilde V_i(x)| \leq \epsilon \obar V_i(x) \qquad \text{and} \qquad |\dtp{c_\perp , \nabla \tilde V_i(x)}| \leq M |\dtp{c_\perp , \nabla \obar V_i(x)}| \,.
\end{equ}
}

\begin{proof}
  We prove the above result in two steps: First of all we show that for certain $\eta_i^*$ the assembled $\tilde V_i$ is a Lyapunov function on $\TT_i$. Then, we construct the sets $\TT_{i}^{(j)}$ such that, by choosing $n_i$ large enough the bounds in \eref{e:disentangle} hold.

  Note that the functions $\obar V_i^{(j)}$ are local Lyapunov functions by \lref{l:local}. Therefore we only have to show that one can safely assemble the local terms $V_i^{(j)}$ to give a Lyapunov function $\tilde V$ on $\TT_i$.
In particular, we proceed to show that for every $\epsilon^* > 0$ there exists a $\eta_i^*> 1$ such that for all $x \in \TT_{i}^{(j)}(c^*)$ and all $r \in \R_c(x)$ we have
\begin{equ}\label{e:ffbound}
  m \Lambda_r(x)|\FF_V^r(x)| \leq \epsilon^*|h_i^{(j)}|\,.
\end{equ}
  By definition of $\obar V_i^{(j)}$ and $\obar V_i^{(j+1)}$ we have that
  \begin{equ}\label{e:boundsonnabla}
    \pd{\nabla_{x,r} \obar V_i^{(j)} - \nabla_{x,r} \obar V_i^{(j+1)}  }  \leq (\eta_i^*-1) \int_{\pi_{ij} (x+ c^r)}^{x+ c^r } \pd{\frac{h_i^{(j)}(z) }{\lambda_{r^*}^w(z)} }\dd z \leq (\eta_i^*-1)c^*\sup_{y \in \BB_{c^*}(z) }\pc{\frac{h_i^{(j+1)}(y)}{\lambda_{r^*}^w(y)}} \,.
  \end{equ}
where we have used the formulation of \eref{e:prop} in the deterministic regime as done in Lemmas~\ref{l:v3}--\ref{l:v1}. Furthermore, using \eref{e:lambdalbub} and the homogeneous scaling behavior of $h_i^{(0)}(x)$ we have that for all $r \in \R$ and $z \in \TT_{ij}(c^*)$ with $\|z\|_2$ large enough
  \begin{equ} \sup_{y \in \BB_{c^*}(z) } \Lambda_r(z)/ \lambda_{r^*}(y) \leq 1 \qquad \sup_{y \in \BB_{c^*}(z) } h_i^{(j+1)}(y)/ h_i^{(j+1)}(z) \leq 2\,. \end{equ}
Combining this result with \eref{e:boundsonnabla} we obtain that there exists a $C>0$ such that
$$\Lambda_r(x)|\FF_V^r(x)| \leq  \Lambda_r (x) \pd{ \nabla_{x,r} \obar V_i^{(j)}-  \nabla_{x,r} \obar V_i^{(j+1)} } \leq C \pc{\eta_i^*-1} h_i^{(j+1)}(x)\,.$$
  Combining this with \lref{l:local} and choosing $\eta_i^* $ close enough to $1$ then establishes \eref{e:ffbound}. Combining this result with \lref{l:local} proves that \eref{e:local} holds on $\TT_i$ for a new constant $C = C-\epsilon^*>0$.

  We now proceed to the second part of the proof. Here, we restrict our attention to $x \in \TT_{ij}$. We immediately obtain the second inequality in \eref{e:disentangle} by noticing that
  $$\dtp{c_\perp, \nabla \obar V_i} = \frac{h_i(x)}{\lambda_{r^*}(x)} = (\eta_i^*)^{n_i}\frac{h_i^{(n_i)}(x)}{\lambda_{r^*}(x)} = (\eta_i^*)^{n_i} \dtp{ c_\perp , \nabla \tilde V_i}\,,$$
  and upon choosing $n_i = \log_{\eta_i^*} M$. For the first part of  \eref{e:disentangle} we write
  \begin{equ} \label{e:chopping}
    \tilde V_i(x) - \tilde V_i(\pi_{ij} x) = \sum_{i = 0}^{n_i} \int_{\pi_{ij} x}^x \indicator_{z \in \TT_i^{(j)}} \frac{ h_i^{(j)}(z) \dd z}{\lambda_{r^*}(x)}  \dd z \leq  \int_{\pi_{ij} x}^x \frac{ h_2^{(0)}(z)\dd z}{\lambda_{r^*}(x)} +   \int_{\pi_{ij} x}^x \indicator_{z \in (\TT_i^{(0)})^c} \frac{h_i^{(0)}(z) \dd z}{\lambda_{r^*}(z)}  \sum_{j = 1}^{n_i}  ((\eta_i^*)^{j}-1)  \,.
  \end{equ}
  Because $\obar V_i(\pi_{ij}x) = \tilde V_i (\pi_{ij}x)$ we obtain the desired result by bounding the second summand on the right hand side of \eref{e:chopping}. In doing so, upon possibly decreasing $\eta_i^*$ further to have $\eta_i^* \in (1,2)$, we write
  \begin{equ}
    |\obar V_i(x) - \tilde V_i (x)| = \sum_{j = 1}^{n_i}  ((\eta_i^*)^{j}-1)  \int_{\pi_{ij} x}^x \indicator_{z \in (\TT_i^{(0)})^c} \frac{h_i^{(0)}(z) \dd z}{\lambda_{r^*}(z)} \leq n_i 2^{n_i} \sup_{y \in \gamma_i(x)} \frac{h_i^{(0)}(y)}{\lambda_{r^*}(y)} \epsilon_i(\TT_i^{(0)},x),,
  \end{equ}
  where $\gamma_i(x)$ denotes the characteristic curve connecting the origin of such curve on the boundary $\TT_{ij'} \neq \TT_{ij}$ to the point $x \in \TT_{ij}$ and we have written $\epsilon_i(\TT_i^{(0)},x) := \int_{\pi_{ij} x}^x \indicator_{z \in (\TT_i^{(0)})^c}\dd z$. We note that $\epsilon_i(\TT_i^{(0)},x) \geq 0$ is decreasing in the size of $\TT_i^{(0)}$, \ie with respect to the partial order induced by the operation of set inclusion, with lower bound $\epsilon_i(\TT_i,x) = 0$.
  Now, by the homogeneous scaling assumptions \eref{e:homogeneous} and \eref{e:homogeneousop}, by the invariance of $\TT_i^{(j)}$ under $\slw$, for $x \in \TT_{ij}$ and because $\slw(\gamma_i(x)) = \gamma_i(\slw(x))$ we have that
\begin{equ}
   {|\obar V_i(\slw(x)) - \tilde V_i (\slw(x))|} \leq  \epsilon_i  l^{\delta_i' - \dtp{c_\inn^{r^*}, w} + 1 } \Big(\sup_{y \in \gamma_i(x)} \frac{h_i^{(0)}(y)}{\lambda_{r^*}(y)} \Big)\,.
\end{equ}
Recall that along $\TT_{ij}$, $\tilde V_i(x)$ and $\obar V_i(x)$ can only differ by a multiplicative constant, so we have
\begin{equ}
   \obar V_i \circ \slw (x) = l^{\delta_i' - \dtp{c_\inn^{r^*}, w} + 1} V_i(x)~.
\end{equ}
Combining this with the boundedness of both $\sup_{y \in \gamma_i(x)} \frac{h_i^{(0)}(y)}{\lambda_{r^*}(y)} $ and $V_i(x)$ we obtain the desired result by choosing $\TT_i^{(0)}$ large enough for $\epsilon(\TT_i^{(0)},x) \leq V_i(x)/\sup_{y \in \gamma_i(x)} \frac{h_i^{(0)}(y)}{\lambda_{r^*}(y)} $.
\end{proof}

\subsection{The assembly process}

We now apply the results obtained in the previous section to show that the local Lyapunov functions from Section\,\ref{s:Lyapunov} can be intuitively assembled to give a global Lyapunov function $V$ on the whole phase space.

\subsubsection{Bulk interfaces: $\TT_{12}$ and $\TT_{23}$}

In this section we assemble the local Lyapunov functions at the interfaces $\TT_{12}$ and $\TT_{23}$. In both cases we can approximate the generator by its dominant transport part $T_i$ defined in \eref{e:transportgen} by \lref{l:taylor2}.

We start by $\TT_{12}$, where by \lref{l:curvature} we can assemble the local Lyapunov functions naturally if we can find parameters for which $\kappa(\obar V, x) < 0$ on $\TT_{12}$. In turn, this condition holds if we can show that
\begin{equ}
\dtp{  c^{r_3} , \nabla \obar V_{1}(x)} > \dtp{c^{r_3}, \nabla \obar V_{2}(x)}\,.
\end{equ}
Writing the terms above as in \eref{e:crnablat12} and using \eref{e:hom2} we obtain the desired inequality through the construction of \lref{l:tuning} for a fixed $\epsilon > 0$ and setting $M > 0$ such that
\begin{equ}
   h_1\pc{5\frac{\delta_1'-5}{b_1} + {1-\delta_1''}} b_1^{\delta_1' - 5}  < M \|c^{r_3}\|_1 h_2 (5b_1+1)^{\delta_2'}b_1^{-5}~.
\end{equ}

Similarly for $\TT_{23}$, we study the convexity of the assembled candidate Lyapunov function $V$ obtained by combining $V_2$ and $V_3$ on the respective regions by considering the $r_3$ directional derivative across the boundary. For $V_2$ and $x \in \TT_{12}$, we have similarly to \eref{e:counter1} that
\begin{equ}
  \dtp{c^{r_3},\nabla \obar V_{2}(x)} = -\|c^{r_3}\|_2 h_2 (5+b_1)^{\delta_2'}b_1^{-2} x_1^{\delta_2'-7}\,.
\end{equ}
Comparing this to the corresponding expression for $\obar V_3$, \ie
\begin{equ}
  \dtp{c^{r_3}, \nabla \obar V_{3}(x) }= \dtp{ c^{r_3},\nabla \pq{h_3 x_1^{\delta_3'-4}x_2^{\delta_3''-2} }_{x_2 = b_1 x_1}} = - h_3(b_2)\pc{5\pc{\delta_1'-4} - \frac{\delta_3''-2}{b_1}} b_1^{\delta_3''-2}x_1^{\delta_3' + \delta_3'' - 7} \,.
\end{equ}
we see that by \eref{e:hom3} we have $\kappa(\obar V, x) < 0$ for $b_1$ large enough if
$$ h_3(b_2)5\pc{\delta_1'-4}b_1^{\delta_3''} > \|c^{r_3}\|_1 h_2 (5+b_1)^{\delta_2'}\,.$$
This in turn holds upon choosing $h_2(b_1,b_2) > 0$ small enough, and we obtain by \lref{l:curvature} that \cndref{a:tanaka} (a) is satisfied at this interface.

 \subsubsection{Boundary interface: $\TT_{00}$, $\TT_{01}$ and $\TT_{34}$}

 In this section we prove that the discrete interface curvature is negative for the interfaces where the continuum approximation \lref{l:taylor2} is not applicable, \ie $\TT_{00}$, $\TT_{01}$ and $\TT_{34}$.

 We start by $\TT_{34}$. In this case, for $x = (b_2,x_2)$ with $b_2$ large enough the discrete curvature term reads
 \begin{equ}
   \obar\kappa_\alpha(\obar V, x) = - h_4 x_2^{\delta_4 } \sum_{k = b_2+1}^{{b_2 + \alpha}}\frac{k^{\delta_4'}}{\delta_4 + \binom{k}{5}5!}   +  h_3 (b_2 + \alpha )^{\delta_3'-4} x_2^{\delta_3''-2} \leq  x_2^{\delta_4}\pc{ h_3{c^* b_2^{\delta_3'-5}} - \frac{h_4b_2^{\delta_4'}}{\delta_4 + \binom{b_2 + c^*}{5}5!}}~, \label{e:012}
\end{equ}
where without loss of generality we consider $\alpha \in (1, c^*)$, in the last inequality we have expanded $(b_2- \alpha)^{\delta_3'-4}$ in $b_2$ and we have applied \eref{e:hom4} and \eref{e:hom3}.
Recalling from \eref{e:hom3} that $h_3(m_4,b_2) = (\delta_3'-4)  m_4 b_2^{-(\delta_3'-4)}$  and our choice of $\delta_3' = \delta_4'$ from \eref{e:deltas}, while bounding $m_4$ from below for $\delta_4' > 5$ as
$$m_4- m_4^* =  h_4 \sum_{k = 1}^{b_2} \frac{k^{\delta_4'}}{\delta_4 + \binom{k}{5}5!} \geq  h_4 \int_1^{b_2-1}k^{\delta_4'-5} \dd k =  \frac{h_4}{\delta_4'-4} ((b_2-1)^{\delta_4'-4}-1)\,,$$
we obtain that to leading order in $b_2$ the right hand side of \eref{e:012} reads,
\begin{equ}
  \obar \kappa_\alpha(x,V) = x_2^{\delta_4}\pc{c^* \frac{h_4}{\delta_4'-4}b_2^{\delta_4'-4}  b_2^{-1} - h_4 b_2^{\delta_4'-5} + \OO(b_2^{\delta_4'-6})}\,.
  \end{equ}
Therefore, for large enough $b_2$ (possibly increasing it from Section~\ref{s:Lyapunov-2}), the discrete interface curvature is negative upon choosing
$$\delta_4' > c^*(\delta_3'-4)+ 4\,.$$

 We now turn to the interface $\TT_{01} := \{x~:~x_2 = b_0\}$. Here we evaluate for $x \in \TT_{01}$
 \begin{equ}
   \obar \kappa_\alpha(\obar V, x) = \obar V_0'(x + e_2) - \obar V_1(x + e_2) = x_1^{\delta_2}
   \Big(h_0' \sum_{k = b_2}^{b_2 + 1} \binom{k}{2}^{-1} - \int_{b_2}^{b_2+1} y^{\delta_1''-2} \dd y \Big)\,.
 \end{equ}
   which has negative sign by choosing $h_0'(b_2,h_1) > 0$ small enough. Combining this with \lref{l:curvature} guarantees that \cndref{a:tanaka} (a) holds at this interface, ensuring natural assembly.


 To conclude, we establish \eref{e:bvp} on the only interface between two regions that do not share a dominant reaction: $\TT_{00} := \{ x \in \Nn_0^2~:~x_2 = 2\}$.
 Because $V_0'$ is not defined for $x_2 = 1$, in this case we have to estimate the cross-term of the differential $\LL V$ explicitly, as we do below.
Setting $x_2 = 2$ we have that
 \begin{equs}
   \LL V (x) &=   \Lambda_3(x) \Delta_3 V_0' + \Lambda_1(x)\Delta_1 V_0'(x) + \Lambda_2(x) \Delta_2 V_0 (x)\,, \label{e:terms}
 \end{equs}
 and we  compute the three summands separately. For the first one we have, for large $x_1$ and $\delta_1 = \delta_1'-5$ that
 \begin{equs}
   \Lambda_3(x) \Delta_3 V_0' & = \Lambda_3(x)\pc{V_0'(\pi_{12 } (x-5e_1)) - V_0'(\pi_{12} x) +\frac{ h_0'}{2} \pc{(x_1-5)^{\delta_1}
   \sum_{k = b_0}^{3} \binom{k}{2}^{-1} - x_1^{\delta_1}\sum_{k = b_0}^{2} \binom{k}{2}^{-1} } }\\
   &  \leq - 2 x_1^{\delta_1'} \frac{h_0'}{2} \pc{ 1 - 5 \delta_1 x_1^{-1}  \sum_{k = b_0}^{3} \binom{k}{2}^{-1}} = - {h_0'}x_1^{\delta_1+5} C_3(x_1,b_0)\,,\label{e:firstterm}
 \end{equs}
for a $C_3$ with $\lim_{x_1 \to \infty}C_3(x_1,b_0)= 1$, where in the last passage we have expanded the difference terms in $x_1$, used that $V_0'(\pi_{12 }x)$ is increasing in $x_1$ and that $\binom{x_1}{5}5! \leq x_1^{5}$. Similarly to \eref{e:firstterm}, we bound the second term in \eref{e:terms} from above by
\begin{equs}
  \Lambda_1(x) \Delta_1 V_0' & = m_1 b_1^{-(\delta_1''-2)} \pc{  (x_1+1)^{\delta_1} - x_1^{\delta_1}   } - h_0' \pc{(x_1+1)^{\delta_1}
  \sum_{k = b_0}^{3} \binom{k}{2}^{-1} - x_1^{\delta_1}\sum_{k = b_0}^{2} \binom{k}{2}^{-1} }\\
  & \leq - x_1^{\delta_1} h_0'  \pc{ 1 + \delta_1 x_1^{-1} \pc{\frac{m_1 }{b_1^{\delta_1''-2}h_0'} + \sum_{k = b_0}^{3} \binom{k}{2}^{-1}}} = - h_0' x_1^{\delta_1} C_1(x_1,b_0,b_1) \,.\label{e:secondterm}
\end{equs}
for a $C_1(x_1,b_0,b_1)$ with $\lim_{x_1 \to \infty}C_1(x_1,b_0,b_1)= 1$ for $b_0, b_1$ fixed.
For the third term in \eref{e:terms} we have by \eref{e:v0}
\begin{equ}
  \Lambda_2(x) \Delta_2 V_0 = \Lambda_2(x)\pc{ V_0(x_1,x_2-1) - V_0(x)} =    x_1^{\delta_0 + \#}( m_0(1)-m_0(2)) \,, \label{e:thirdterm}
\end{equ}
where $\#$ is $0$ for \abbr{crn}0 and $1$ for \abbr{crn}1. Therefore, recalling that $\delta_0 = \delta_1$ from \eref{e:deltas}, a comparison of the exponents of \eref{e:firstterm}--\eref{e:thirdterm} shows that for large $x_1$ and $x \in \TT_{00}$, we have $\LL V(x) \leq -(1-\epsilon(x_1,b_2)) h_0'(x)$ for $\lim_{x_1 \to \infty}\epsilon(x_1,b_2) = 0 $ as required.

\subsection{Summary}
Given the involved nature of the process carried out in the previous sections, we summarize for the reader the order of choice of parameters for the construction of our global Lyapunov function $V$. First off, one starts by setting all the powers of the homogeneous functions $\{V_i, h_i\}$, \ie $\delta_i, \delta_i'$ (and $\delta_i''$ when available) as in \eref{e:deltas}, except for $\delta_4'$ that will be fixed later. Then one proceeds to check that the propagated Lyapunov functions satisfy \eref{e:bvp} in the regions $\TT_{i}$, thereby fixing large enough constants $\{b_i\}$. Finally, one has ensure that the natural assembly condition \eref{e:tanaka} holds on the interfaces $\TT_{ij}$.
This last step is performed sequentially from the priming region through the transport regions to the diffusive regions and fixes, in the order, $\delta_4'$, $h_4$, $h_3(b_2)$, $h_2(b_1,h_3)$, $h_1(h_2,b_1)$, $h_0'(h_1,b_0)$, $h_0(h_0',b_0)$. Recall that by \eref{ex:1} we also have to construct the regions $\TT_{2}^{(j)}$. During this procedure we choose $M(h_1, b_1)$ large enough and consequently $n_2(M,\eta_2^*)$.
In doing so, the parameter $h_1$ is increased by a constant $\epsilon$ small at will, so this has only marginal effect on the choice of the subsequent parameters. Finally, choosing a large enough $\rho > 0$ ensures that \eref{e:phi} holds for $\|x\|_2> \rho$.

\section{Large time asymptotics}

\label{s:convergence}

This section is dedicated to the study of the invariant measure of the \abbr{crn}s \eref{e:e0} and \eref{e:e1}.
\cf{Based on the structure of the Lyapunov function solving \eref{e:bvp}, we establish asymptotic upper bounds on the density of the invariant measure and estimate the speed of convergence to such invariant measure. To do so we impose some conditions on the stability and mixing properties of the process $X_t$. In particular, in terms of stability we require that \cndref{cnd:stability} holds, \ie that the process returns, on average, to a compact set in $\mathbb N_0^d$ by bounding its average speed towards that set. Similarly to \cite{hairer16, hm11} we then bound the convergence to the invariant measure from above by the probability that two processes with the same generator couple in a ``small'', attracting region of phase space. To obtain this result, we need to ensure that such ``small'' set exists and that it is sufficiently mixing, as encoded in \cndref{cnd:mixing}.}

\subsection{Invariant measure density}

We now proceed to prove the existence and uniqueness of the invariant measure for a process $X_t$ satisfying \cref{cnd:stability} through the following, standard result.
\lmm{[\cite{hm11, meyn12}] \label{l:invmeas} Under Conditions~{cnd:stability} and \ref{cnd:mixing} there exists a $\sigma$-finite invariant measure $\mu$ for the process $X_t$. Furthermore, if $\vphi > 1$ then there exists a constant $\Cmu > 0$ such that
 \begin{equ}
  \int \vphi(V(x)) \mu(\dd x) \leq \Cmu\,.\label{e:invmeas}
 \end{equ}
}
\begin{proof}
 The existence of an invariant measure under the assumptions of the lemma was established in \cite[Theorem 12.3.3]{meyn12}. The finiteness of the integral in \eref{e:invmeas} is established in \cite[Theorem 14.0.1]{meyn12} and in \cite{hm11} for the case of positive recurrent process, as for \abbr{crn0}. We adapt the proof to the null recurrent case below.

 Let $f \in C_b(\RR^2)$ be compactly supported. Then, denoting for any set $B \subset \Nn_0^d$ by $\tau_B$ the first return time to that set and by $m_B$ an probability distribution on $B$ at time $t = 0$ we define  for any $N > 0$ as in \cite[Theorem 3.5.3.]{norris98} the measure $\nu_N$ by
 \begin{equ}\label{e:nuN}
   \nu_N(A) :=  \Ex{m_B}{\int_{0}^{N\wedge \tau_B - 1} \indicator_A(X_s)  \dd s}\,.
 \end{equ}
 Then, using assumption \eref{e:phi} and Dynkin's formula we obtain, for a compact $B \supseteq \KK$,
\begin{equs}
  \nu_N \vphi(V(x)) & = \Ex{m_B}{\int_0^{\tau_B \wedge \tau_{\{\|x\|_2 \geq N\}} } \vphi(V(X_s)) \dd s  } \leq - \Ex{m_B}{\int_0^{\tau_B \wedge \tau_{\{\|x\|_2 \geq N\}} \wedge N } \LL V(X_s) \dd s  }\\& \leq \Ex{m_B}{V(X_0) } \leq \sup_{x \in B} V(x) \leq \Cmu\,.\label{e:finitemeas}
\end{equs}
  The uniform in $N$ finiteness of the upper bound \eref{e:finitemeas} and the convergence of $\nu_N$ to the invariant measure of $X_t$  for $N \to \infty$ proves the desired result by application of Fatou's lemma.
\end{proof}

\lref{l:invmeas} allows to bound the invariant measure density from above: by finiteness of \eref{e:invmeas} and since $\int \|x\|^{-d-\epsilon}\dd x < \infty$ we can bound from above the tails of the invariant measure by writing, for any $\epsilon > 0$,
$$ \mu(\dd x) < \pc{\|x\|^{d+\epsilon} \vphi(V(x)}^{-1} \dd x~.$$

The exponential convergence to the invariant measure in the case of \abbr{crn1} is a standard result \cite{hm09} under Assumptions~\ref{cnd:stability} and \ref{cnd:mixing}, provided that there exist $\gamma > 0$ such that $h(x) > \gamma V(x)$. To see that this condition is satisfied in our case, it is sufficient to note that in Lemmas~\ref{l:v4}--\ref{l:v0} the scaling exponent of $h_i$ is always larger or equal than the one of $V_i$ in all the regions $\TT_i$.

\cf{
\subsection{Convergence to the invariant measure}

We now establish results for the rate of convergence to the invariant measure. This extends known results for concave $\vphi(x)$ \cite{hairer16, meyn12} to the case $\vphi(x) \to 0$. To obtain such results, we compare the process $V(X_t)$ to a deterministic flow $\Phi(u,t)$ solving the system of ordinary differential equations
\begin{equ}
 \partial_t \Phi(u,t) = - \tilde \vphi(\Phi(u,t)) \qquad \text{with} \qquad \Phi(u,0) = u\,,
 \label{e:flow}
\end{equ}
where $\tilde \vphi(x) = \Cphi \vphi(x)$ for an arbitrary constant $\Cphi \in (0,1)$ to be chosen later.
To do so, we introduce the function
\begin{equ}
 \hfu = \int_1^u \frac{\dd s}{\tilde \vphi(s)}~,
 \label{e:Hphi}
\end{equ}
returning the travel time of a particle in the flow \eref{e:flow} to the point $u = 1$.

\aa{We define
 \begin{equ}
  \PP_t^G(x,y) := \frac{\PP_t(x,y)}{\int \PP_t(x,y)G(y)\dd y}\qquad \text{and} \qquad \mu^G(\dd y) := \frac{\mu(\dd y)}{\int G(x) \mu(\dd x)}\,.
 \end{equ}
 We also define the weighted total variation metric $$d_G(\nu_1,\nu_2) := \sup_{\phi:\Rr \to \Rr, |\phi|\leq G}\pq{\int \phi(x) \nu_1(\dd x)- \int \phi(x) \nu_2(\dd x)}\,,$$
 and the corresponding total variation norm $\tvn{\cdot}^G$.}
\begin{pro}
 {Under Conditions~\ref{cnd:stability} and \ref{cnd:mixing} there exists a constant $C>0$ such that for every $x, y \in \mathbb N_0^d$ and all $t > 0$ we have}
 \begin{equ}
  \tvn{\PP_t(x, \,\cdot\,) - \PP_t(y, \,\cdot\,)} \leq C \frac{V(x)+V(y)}{\hf1t}~. \label{e:convergence1}
 \end{equ}
 \aa{Furthermore, assume that there exists a strictly increasing $h(x) \sim \Omega(x^\epsilon)$ and a function $G(x)\colon \RR^d \to \RR$ such that
  \begin{equ}
   \int h(V(x)) G(x) \mu(\dd x) = K' < \infty~.
  \end{equ}
  Then there exists a constant $C' > 0$ such that
  \begin{equ}
   \tvn{\PP_t^G(x,\,\cdot\,) - \mu^G}^G \leq \frac{C V(x)}{\hf1t} + \frac{C'}{(\hf1t)^{\kappa}}\,,
   \label{e:convergence2}
  \end{equ}
  where $\kappa = \frac{1}{1+1/\delta \epsilon}$.}
 \label{p:convergence}
\end{pro}
We prove the above through the following intermediate results.


\lmm{\label{l:lphi} Let $\vphi$ be convex and bounded from above and let $V$ be concave, then there exists $\Cphi' > 0$ such that, for any $\Cphi < 1$ and for all $t > 0$\,,
 \begin{equ}\label{e:lphi}
  \LL \Phi^{-1}(V(X_t),t) \leq \Cphi'\partial_u \Phi^{-1}(V(X_t),t)\, \LL V(X_t)~,
 \end{equ}
 where $\Phi^{-1}(u,t)$ is the inverse in the first coordinate of the flow \eref{e:flow}.}

\proof{Denoting by $D_x$ the total derivative in the argument of $V(\,\cdot\,)$ we rewrite the difference term in the generator as
 \begin{equs}
   \Delta_r \Phi^{-1}(V(X_t),t) =  \Phi^{-1}(V(X_t+c_r),t) - \Phi^{-1}(V(X_t),t) = \int_0^{c^r} D_x \Phi^{-1}(V(X_t+y),t) \dd y \,. \label{e:remainder}
 \end{equs}
  Proceeding to analyze the right hand side of the above expression, we note that we can write
\begin{equ}
  \Phi^{-1}(u,t) = H_\vphi^{-1}(H_\vphi(u)+t)~.\label{e:h-1}
\end{equ}
Combining \eref{e:h-1} with the fact that $\partial_u H_\vphi = 1/\tilde \vphi$ and \eref{e:flow} we obtain
\begin{equ}
  \partial_u\Phi^{-1}(u,t) = \frac{\partial_t \Phi^{-1}(u,t)}{\tilde \vphi(u)} = \frac{\tilde \vphi( \Phi^{-1}(u,t))}{\tilde \vphi(u)}\,.
  \end{equ}
  In light of this, there exist constants $C, C' > 0$ such that
\begin{equs}
   \pd{\frac{\sup_{y \in \BB_{c^r}(x)}D_x (\partial_u \Phi^{-1})(V(y),t)}{\partial_u \Phi^{-1}(V(x),t)}}
  & = \sup_{y \in \BB_{c^r}(x)} \partial_x V(y) \frac{\tilde \vphi'(\Phi^{-1}(V(y),t)) - \tilde \vphi'(V(y))}{\tilde \vphi(V(y))} \frac{ \partial_u \Phi^{-1}(V(y),t)}{\partial_u \Phi^{-1}(V(x),t)}
  \\&\leq C \sup_{y \in \BB_{c^r}(x)} \partial_1 V(x) \frac{ \vphi'(\Phi^{-1}(V(y),t)) -  \vphi'(V(y))}{ \vphi(V(y))} \leq C'\,,\label{e:remainder2}
\end{equs}
where in the first inequality we have used that $\Phi$, like $\vphi$, is convex and bounded, and in the second that $V$ is concave, $\Phi$ is increasing in $t$ and $\vphi'/\vphi \to 0$\,.
Combining \eref{e:remainder} with \eref{e:remainder2} we obtain
\begin{equ}
\Delta_r \Phi^{-1}(V(X_t),t) \leq \pc{1+R\pc{c^r,V(X_t)}} \partial_u \Phi^{-1}(V(X_t),t) \Delta_r V(X_t)~,
\end{equ}
with a remainder $R(c^r,V(X_t))$ that is bounded from above by Taylor's theorem. This concludes the proof by choosing $\Cphi' := (1+ \sup_{x \in \RR^2} R(c^r,V(X_t)))$.\qed
}

\lmm{\label{l:hitting}For any $K > 0$, define $\tau_K := \inf\{t > 0~:~ V(X_t) \leq K\}$. Then there exists constants $\Cphi \in (0,1) $ and $\Ctau > 0$ so that for all $x \in \RR^2$ we have
 \begin{equ}\label{e:hitting}
  \Ex{x}{\hft{\tau_K}} \leq \Ctau V(x)~.
 \end{equ}
}

\proof{Consider the process $F(X_t) := \Phi^{-1}(V(X_t),t)$. By Dynkin's formula we have
 \begin{equs}
  \Ex{x}{\hft{t}} & = V(x) + \Ex{x}{\int_0^t \pc{\LL \hft{s} + \partial_s \hft{s}} \dd s}\,,
 \end{equs}
 and proceed to show that the term in the integral is negative. This results by application of \lref{l:lphi} to the differential of $\Phi^{-1}(V(X_t),t)$:
 \begin{equs}
  \LL \hft{s} + \partial_s \hft{s}&  \leq \Cphi' \partial_u \Phi^{-1}(V(X_t),t) \LL V(X_t) + \partial_t \Phi^{-1}(V(X_t),t)\\
  & \leq -\Cphi' \frac{\partial_t \Phi^{-1}(V(X_t),t)}{\Cphi \vphi(V(X_t))} \vphi(V(X_t)) + \partial_t \Phi^{-1}(V(X_t),t) = 0\,,
 \end{equs}
 upon choosing $\Cphi = \Cphi'$ small enough.
 \qed
}

\lmm{\label{l:coupling}
 Let $X_t, Y_t$ the two processes described above and let $\tau_C$ be their coupling time, \ie $\tau_C := \inf\{T > 0 ~:~X_t = Y_t\, \forall t>T\}$. Then we have that
 \begin{equ}\label{e:coupling}
  \tvn{\PP_t(x,\,\cdot\,) - \PP_t(y,\,\cdot\,)} \leq C\frac{V(x)}{\hf1t}~.
 \end{equ}
}

\begin{proof}{To obtain the above result, we separate the events $\{\tau_C > \tau_n\}$ and $\{\tau_n>t\}$
  \begin{equs}
   \Ex{x}{\phi(X_t)} & = \Ex{x}{\phi(X_t) \pc{I_{\{\tau_n \leq t\}} + I_{\{\tau_n>t\}}(I_{\{\tau_C > \tau_n\}} + I_{\{\tau_C \leq \tau_n\}}) }}\\
   & \leq \Ex{x}{\phi(Y_t) I_{\{\tau_n \leq t\} \cup \{\tau_C \leq \tau_n\}}} + \Ex{x}{\phi(X_t) I_{\{\tau_n>t\}}} + \Ex{x}{\phi(X_t) I_{\{\tau_C > \tau_n\}}}~.
  \end{equs}
  The probability of the above events can be bounded from above by $\p{\tau_C > \tau_n} \leq (1-\alpha)^n $ and
  \begin{equ}
   \px{x}{\tau_n>t} \leq \px{x}{\tau_1>t/2} + \px{V_C}{\tau_n-\tau_{1}> t/2} \leq \frac{V(x)}{\hf1t} + 2(n-1)\frac{\Ex{V_C}{\Delta \tau}}{t}~.
  \end{equ}
  Choosing $n = \frac{ \log {V(x)}/{\hf1t}}{\log(1-\alpha)}$ we can write
  \begin{equs}
   \tvn{\PP_t(x,\,\cdot\,) - \PP_t(y,\,\cdot\,)} & = \sup_{\phi \in L^\infty} \pq{\int \phi(z) \PP_t(x,\dd z) - \int \phi(y) \PP_t(y,\dd z)} = \sup_\phi \pq{ \Ex{x}{\phi(X_t)} - \Ex{y}{\phi(Y_t)}}\\
   &\leq \|\phi\|_\infty \pq{\px{x}{\tau_n>t} + \px{x}{\tau_C > \tau_n}} \leq C \frac{V(x)}{\hf1t}~.
  \end{equs}
 }
\end{proof}

\aa{\lmm{\label{l:convergence}Let $h(x) \sim x^{-\epsilon}$ for $\epsilon > 0$\,. Then \eref{e:convergence2} holds with $\kappa = \pc{1+1/\delta\epsilon}^{-1}$. }

 \begin{proof}
  Using that
  \begin{equs}
   G(y) \mu^G(\dd y) = \int \PP_t(z, \dd y) G(y) \mu^G(\dd z) = {\int   {\int G(w)\PP_t(z, \dd w) } } G(y) \PPG(z, \dd y)  \mu^G(\dd z)\,,
  \end{equs}
  we separate the case $z > y$ and the complementary case: by the decreasing property of $G(y)$ we have
  \begin{equs}
   \int G(y) \pd{\mu^G(\dd y) - \PPG(x,\dd y)} = &  \int \mu^G(\dd z)  {\int_{V > R} G(w)\PP_t(z, \dd w) }\int G(y)  \pd{\PPG(z, \dd y) - \PPG(x,\dd y)}\\
   &+  \int \mu^G(\dd z) {\int_{V < R} G(w)\PP_t(z, \dd w) }\int G(y) \pd{\PPG(z, \dd y) - \PPG(x,\dd y)}\,,
  \end{equs}
  We then consider the three summands separately. For the first term we write, by triangle and Chebyshev inequality,
  \begin{equ}
   \int\mu^G(\dd z)  {\int_{V > R} G(w)\PP_t(z, \dd w) } \int G(y)\pd{\PPG(z, \dd y) - \PPG(x,\dd y)} \leq 2 \mu^G(V > R) \leq 2 \frac{\int G(y)h(y)\mu^G(\dd y)}{h(R)}~.
  \end{equ}

  To bound the second term from above, we split the integral over $\dd z$ in the regime where $V(z)\geq R$ and its complement. Respectively, we obtain
  \begin{equs}
   \int_{V > R} \mu^G(\dd z) {\int_{V<R} G(w)\PP_t(z, \dd w) }\int G(y) \pd{\PPG(z, \dd y) - \PPG(x,\dd y)} \leq\,,
  \end{equs}
  and
  \begin{equs}
   \int_{V < R} \mu^G(\dd z) {\int_{V<R} G(w)\PP_t(z, \dd w) }\int G(y) \pd{\PPG(z, \dd y) - \PPG(x,\dd y)}\,,
  \end{equs}

  Then we obtain, by Fubini-Tonelli and triangle inequality,
  \begin{equs}
   {\int \pd{\PPG(z, \dd y) - \mu^G(\dd y)} G(y)} &= \int  \mu^G(\dd z) \int \PP_t(z, \dd w) G(w) \int G(y)\pd{\PPG(x,\dd y) - \PPG(z, \dd y)}\\
   & = \int \mu^G(\dd z) \int G(y) \pd{\PP_t(x,\dd y) - \PP_t(z, \dd y) + \PP_t(x, \dd y) \pq{\frac{\int G(w) \PP_t(z, \dd w)}{\int G(w) \PP_t(x, \dd w)}-1}}\\
   & \leq 2 \int  \mu^G(\dd z)  \int G(y)\pd{\PP_t(x,\dd y) - \PP_t(z, \dd y)}~.
  \end{equs}
  Consequently, for any $R> 0$,we have
  \begin{equs}
   \tvn{\PPG(x, \,\cdot\,) - \mu^G}^G &
   \leq \int \mu^G(\dd z)  \pq{ \int_{V < R} G(y)\pd{\PP_t(x, \dd y) - \PP_t(z, \dd y)}+\int_{V > R} G(y)\pd{\PP_t(x, \dd y) - \PP_t(z, \dd y)}}\\
   & \leq \obar G \int \mu^G(\dd z) \tvn{\PP_t(x, \,\cdot\,) - \PP_t(z, \,\cdot\,)} + 2\int_{z > x} \int_{V> R} G(y) \mu^G(\dd y) + \int_{z < x}\int_{V> R} G(y) \mu^G(\dd y)\\
   & \leq  \obar G \int \mu^G(\dd z) \tvn{\PP_t(x, \,\cdot\,) - \PP_t(z, \,\cdot\,)} + 2\int_{V> R} G(y) \mu^G(\dd y)~,
  \end{equs}
  where in the last step we have used that $|\PP_t(x, \dd y) - \PP_t(z, \dd y)| \leq 2\PP_t(z,\dd y)$ for large enough $z$.\\
  Now, consider the two summands separately. For the first
  \begin{equs}
   \int_{V< R} \pd{\PP_t^G(x,y)-\mu^G(y)} \dd y \leq \int_{V< R} \tvn{\PP_t(x,\,\cdot\, )-\PP_t(y,\,\cdot\,)} \mu^G(\dd y) \leq \frac{C}{\hf1t} \pc{V(x) + \int_{V < R} V(y) \mu(\dd y)}\,,
  \end{equs}
  which can be bounded from above by noticing that $x/\vphi(x)$ is an increasing function, so $V(y) \leq \vphi(V(y)) R/\vphi(R)$ and since $\int \vphi(V(x)) \mu(\dd x) < C$ we have
  \begin{equ}
   \int_{V< R} \pd{\PP_t^G(x,y)-\mu^G(y)} \dd y \leq \frac{C}{\hf1t} \pc{V(x) + \frac{C R}{\vphi(R)}}~.
  \end{equ}
  Considering now the second term, for any strictly increasing function $h: \RR \to \RR$ growing as $x^\epsilon$ we have, by Markov inequality and by our assumption on $G$
  \begin{equ}
   \int_{V> R} G(y) \mu^G(\dd y) =  \int_{h(V)> h(R)} G(y) \mu^G(\dd y) \leq h(R)^{-1} \int h(y) G(y) \mu^G(\dd y) \leq C h(R)^{-1}\,.
  \end{equ}
  Combining (a) and (b) we have that
  \begin{equ}
   \tvn{\PP_t^G-\mu^G}^G \leq \frac{C}{\hf1t} \pc{V(x) + \frac{C R}{\vphi(R)}} + C h(R)^{-1}~.
  \end{equ}
  Optimizing over $R$ yields $R \sim \pc{\hf1t}^{\frac{1}{\epsilon + 1/\delta}}$, which if inserted in the above expression yields
  \begin{equ}
   \tvn{\PP_t^G-\mu^G}^G \leq \frac{C}{\hf1t} V(x) + \frac{C'}{\pc{\hf1t}^{\frac{1}{1 + 1/\delta\epsilon}}} ~.
  \end{equ}
\end{proof}}

\aa{left to do:\\
 verify the above proof\\
 adjust the proof with end in red above on the indicators\\
 add the dual $W(x,y)$ using that $$\phi(V(x))+\phi(V(y)) \leq \phi(V(x) + V(y)) $$}

\begin{proof}[Proof of \pref{p:convergence}]
 \eref{e:convergence1} follows by combining $\tvn{\PP_t(x,\,\cdot\,)-\PP_t(y,\,\cdot\,)} \leq \px{(x,y)}{\tau_C>t}$ with \lref{l:coupling}
\end{proof}}

\section*{Acknowledgements}
AA acknowledges the support of the Swiss National Science Foundation through the grant P2GEP2-175015. JCM acknowledges the support of the U.S. National Science Foundation  through the grant DMS-1613337. Furthermore, AA would like to thank Noé Cuneo for interesting discussions.

\bibliographystyle{JPE.bst}
\bibliography{bib.bib}

\appendix

\section*{Appendix}

\section{Explitic estimates}

\begin{proof}[Diffusive Region]
In this section we evaluate \eref{e:prop} in the region $\TT_0$. We do so by considering the sets $\{x_1 = 0\}$ and $\{x_1 = 1\}$ separately.

 For both \abbr{crn0} and \abbr{crn1} we see that for $x = (x_1,0)$ the first term of \eref{e:prop} reads
 \begin{equ}
   \Ex{(x_1,0)}{V_0'(X_\tau)} = \Ex{(x_1+1,1)}{V_0'(X_\tau)}\,,
   \end{equ}
   while for the second term we have
   \begin{equ}
     \Ex{(x_1,0)}{\int_0^{\tau} h_0(X_t)\, \dd t} = h_0 x_1^{\delta_0'}\E{\Delta \tau_{(x_1,0)}} + \Ex{(x_1+1,1)}{\int_0^{\tau} h_0(X_t)\, \dd t}\,.
   \end{equ}
  Therefore, we obtain
  \begin{equ} \label{e:level0}
    V_0((x_1,0)) = V_0(x_1+1,1) + h_0 x_1^{\delta_0'}\,,
  \end{equ}
thereby reducing the problem to the computation of $V_0$ on the set $\{x_2 = 1\}$, as we do in the next paragraph.

For \abbr{crn0} we calculate the first expectation in \eref{e:prop} for $x = (x_1,1)$ by writing
 \begin{equ}
  \Ex{(x_1,1)}{m_0(X_\tau)_1^{\delta_0}} = \sum_{k = x_1}^\infty \Ex{(x_1,1)}{m_0(X_\tau)_1^{\delta_0} | \uparrow_k} \px{(x_1,1)}{\uparrow_k} =  \sum_{k = x_1}^\infty k^{\delta_0} (1-a) a^{k-x_1}\,,\end{equ}
where $a = \p{\downarrow_k} < 1$ and the definition of the event $\downarrow_k$ was given in Claim~\ref{cla:transiant}.
  Bounding the right hand side from above and below by integration produces
  \begin{equ}
     \int_1^\infty e^{k \log a} (k + x_1)^{\delta_0} \dd k \leq \frac{\Ex{x}{m_0(X_\tau)_1^{\delta_0}}}{m_0 (1-a)} \leq \int_0^\infty e^{k \log a} (k + x_1)^{\delta_0} \dd k \,.
  \end{equ}
From this we see that for \abbr{crn0} the first term in \eref{e:prop} is well-defined for all choices of $\delta_0$ and behaves asymptotically as $V_0 \ssim{e_1} l^{\delta_0}$. We will keep only the term corresponding to this scaling for the formula \eref{e:v0}, as the remaining terms are negligible in the scaling of interest.

  Similarly, for \abbr{crn1} we have by \eref{e:pupcrn1}
  \begin{equ}
   \Ex{(x_1,1)}{m_0(X_\tau)_1^{\delta_0}} = m_0\sum_{k = x_1}^\infty \Ex{(x_1,1)}{(X_\tau)_1^{\delta_0} | \uparrow_k} \px{(x_1,1)}{\uparrow_k} =  \sum_{k = x_1}^\infty k^{\delta_0} \frac{x_1}{k(k+1)} ~.\end{equ}
   The sum on the right hand side of the above expression can be bounded from above and below by integration
   \begin{equ}
        x_1 \int_1^\infty (k + x_1)^{\delta_0-2} \dd k  \leq \Ex{(x_1,1)}{(X_\tau)_1^{\delta_0}} \leq x_1 \int_0^\infty (k + x_1)^{\delta_0-2} \dd k~.
      \end{equ}
    By the divergence of the integral on the left hand side for $\delta_0 \geq 1$  and the convergence of the one on the right hand side for $\delta_0 < 1$ we obtain the desired well-definiteness result for $V_0$. The dominant scaling behavior for this component of the sum is given by $V_0 \ssim{e_1} l^{\delta_0}$ by the result of the two integrations above.

   We now proceed to estimate the second term on the right hand side of \eref{e:prop}. Using the scaling properties of $h_0(\,\cdot\,)$ and using the linearity of the expectation we have
   \begin{equs}
    \Ex{(x_1,1)}{\int_0^\tau h_0 \pc{X_t} \dt} &= \sum_{k = x_1}^\infty \Ex{(x_1,1)}{\int_0^\tau h_0 \pc{X_t} \dt | \uparrow_k} \px{(x_1,1)}{\uparrow_k}\\&
    = \sum_{k = x_1}^\infty \Ex{x}{\sum_{l = x_1}^{k} h((l,1)) \Delta \tau_{(l,1)} + h((l,0)) \Delta \tau_{(l,0)} | \uparrow_k} \px{(x_1,1)}{\uparrow_k}  \\&
    =  \sum_{k = x_1}^\infty \sum_{l = x_1}^k l^{\delta_0'} \E{\Delta \tau_{(l,1)} + \Delta \tau_{(l,0)}} \px{(x_1,1)}{\uparrow_k}~. \label{e:eintxdelta}\end{equs}

    Evaluating \eref{e:eintxdelta} for \abbr{crn0} with $h_0(x) = h_0 x_1^\dpr$ yields
    \begin{equs}
     \Ex{(x_1,1)}{h_0\int_0^\tau (X_t)_1^\dpr \dt} &
     = h_0\sum_{k = x_1}^\infty \sum_{l = x_1}^k l^\dpr \E{\Delta \tau_{(l,1)} + \Delta \tau_{(l,0)}} \px{(x_1,1)}{\uparrow_k}\\
     & = h_0 a^{-x-1} \sum_{k = x_1}^\infty e^{-\log(1/a) k} \sum_{l = x_1}^k l^\dpr \pc{\frac{1}{\kappa_1 + \kappa_2} + \frac{1}{\kappa_1}}\Bigg|_{\kappa_1=\kappa_2=1}\\
     & = \frac{3}{2}h_0 a^{-x_1-1} \sum_{k = x_1}^\infty e^{-\log(1/a) k} \sum_{l = x_1}^k l^\dpr\,,
    \end{equs}
    Now, similarly to what was done for the first term of \eref{e:prop} we bound the right hand side of the above equation from above and below by integration:
    \begin{equ}
      \int_1^\infty \dd k\, e^{-\log(1/a) k} \int_0^{k} \dd l\, (x_1+l)^\dpr \leq  \frac{\Ex{(x_1,1)}{h_0\int_0^\tau (X_t)_1^\dpr \dt}}{h_0}\leq \int_0^\infty \dd k\, e^{-\log(1/a) k}\int_0^{k} \dd l\, (x_1+l)^\dpr\,.
    \end{equ}
    Bounding the right hand side from above by $x^\dpr + \int_0^\infty \dd k\, e^{-k} k^{\dpr + 1}$ we see that it is well defined for all values of $\dpr \in \Rr$ and it scales as $l^{\dpr}$ to leading order.

    Similarly, for \abbr{crn1} we rewrite \eref{e:eintxdelta} using \eref{e:pupcrn1} as
    \begin{equs}
     \Ex{(x_1,1)}{h_0\int_0^\tau \pc{X_t}_1^\dpr \dt} &
     = h_0 \sum_{k = x_1}^\infty \sum_{l = x_1}^k l^\dpr \Ex{l}{\Delta \tau_{(l,1)} + \Delta \tau_{(l,0)}} \px{x_1}{\uparrow_k}\\
     & = h_0 \sum_{k = 0}^\infty \frac{x_1}{(x_1+k)(x_1+k+1)} \sum_{l = x_1}^{k} l^\dpr \pc{\frac{1}{l+1} + 1}\,.
    \end{equs}
    We bound from above and below (up to an approximating constant) the right hand side by
    \begin{equ}
      x_1 \int_0^\infty \dd k\, \frac{1}{(x_1+k)^2} \int_{x_1}^{x_1+k} \dd l\,l^\dpr = \frac{x_1}{\delta_0'} \int_{0}^\infty \dd k\, \frac{x_1^{\dpr+1} - (k+x_1)^{\dpr + 1}}{(k+x_1)^2} \,. \label{e:level11}
    \end{equ}
    We notice that the right hand side of the above expression converges only if $\dpr < 0$. If this condition holds we have that $V_0 \ssim{e_1}l^{\delta_0'+1}$.

    Combining \eref{e:level0} and the scaling behavior of \eref{e:prop} on $\{x~:~(x_1,1)\}$ results in \eref{e:v0}.

    \end{proof}

    \begin{proof}[Priming region]
     As in the case of the diffusive region, we couple the process at hand with a new process, denoted by $ X_n$, jumping to the left ($\leftarrow$) or down ($\downarrow$) with the respective site-dependent probabilities
     $$\px{x}{\leftarrow_{x}} = \frac{x_2^2 \binom{x_1}{5}5!}{x_2^3 + x_2^2 \binom{x_1}{5}5!} = \frac{ \binom{x_1}{5}5!}{x_2 + \binom{x_1}{5}5!} \quad \text{and} \quad \px{x}{\downarrow_{x}} =\frac{x_2^3}{x_2^3 + x_2^2 \binom{x_1}{5}5!} = \frac{x_2 }{x_2 + \binom{x_1}{5}5!}\,. $$
      Consequently, the hitting distribution of the set $ \{x~:~x_1 = b_2-1\}$ by the process $X_n$ with initial condition $X_0 = (b_2,x_2)$ and $x_2'\leq x_2$ is given by:
     \begin{equs}
      \px{(b_2,x_2)}{\leftarrow_{(b_2,x_2')}} &= \px{(b_2,x_2)}{\leftarrow_{(b_2,x_2')}|\downarrow_{(b_2,x_2 \dots x_2')}} \px{(b_2,x_2)}{\downarrow_{(b_2,x_2 \dots x_2')}} \\&= \px{(b_2,x_2')}{\leftarrow_{(b_2,x_2')}} \prod_{k = x_2}^{x_2'-1} \px{(b_2,k)}{\downarrow_{(b_2,k)}} = \\
      & = \frac{\binom{b_2}{5}5! }{x_2' + \binom{b_2}{5}5!}\prod_{k = x_2'+1}^{x_2} \frac{k }{k + \binom{b_2}{5}} = \frac{B_5(b_2)}{x_2'+B_5(b_2)} \frac{\Gamma(x_2+1)}{\Gamma(x_2'+1)}\pc{\frac{\Gamma{(x_2+B_5(b_2)+1)}}{\Gamma(x_2'+B_5(b_2)+1)}}^{-1} \\
      & = B_5(b_2)\frac{\Gamma(x_2'+B_5(b_2))}{\Gamma(x_2'+1)}\frac{\Gamma(x_2+1)}{\Gamma{(x_2+B_5(b_2)+1)}}\,,
     \end{equs}
     where $\Gamma(\,\cdot\,)$ is the gamma function and we have written $B_5(x) := \binom{x_1}{5}5!$. For large values of $x_2, x_2'$ we approximate
     \begin{equ}\label{e:approxproba}
       \px{(b_2,x_2)}{\leftarrow_{(b_2,x_2')}} \sssim{e_2} B_5(b_2) x_2'^{B_5(b_2)-1} x_2^{-B_5(b_2)}\,,
       \end{equ}
       where from now on the symbol $\sssim{e_2}$ denotes equality up to subdominant terms in $\slw$ for $w = e_2$.
       Iterating over multiple levels from $b_2$ to $0$, we obtain for the first term of \eref{e:prop} for $V_4(\, \cdot \,)$ scaling homogeneously as $V_4 \ssim{e_2} l^{\delta_4}$ by
     \begin{equs}
      \frac{\Ex{(b_2,x_2)}{V_4(X_\tau)}}{m_4^*} &= \sum_{x_2' = 0}^{x_2} \px{(b_2,x_2)}{X_\tau = x_2'} (X_\tau)_2^{\delta_4} \\
      &= \sum_{k_{b_2} = 0}^{x_2} \px{{(b_2,x_2)}}{\leftarrow_{(b_2,k_{b_2})}} \sum_{k_{b_2-1} = 0}^{k_{b_2}} \px{{(b_2-1,k_{b_2})}}{\leftarrow_{(b_2-1,k_{b_2-1})}} \cdots \sum_{k_1 = 0}^{k_2} \px{{(1,k_2)}}{\leftarrow_{(1,k_{1})}} k_1^{\delta_4}
      \\&= \sum_{k_{b_2} = 0}^{x_2} B_5(b_2)\frac{\Gamma(x_2+1)}{\Gamma{(x_2+B_5(b_2)+1)}}\frac{\Gamma(k_{b_2}+B_5(b_2))}{\Gamma(k_{b_2}+1)} \sum_{k_{b_2-1} = 0}^{k_{b_2}} B_5(b_2-1)\frac{\Gamma(k_{b_2}+1)}{\Gamma{(k_{b_2}+B_5(b_2-1)+1)}}\\&\qquad \qquad\cdots \sum_{k_1 = 0}^{k_2} B_5(1) \frac{\Gamma(k_{2}+1)}{\Gamma(k_{2}+B_5(1)+1)}\frac{\Gamma(k_{1}+B_5(1))}{\Gamma(k_{1}+1)}  k_1^{\delta_4}\\
      &=\prod_{j = 1}^{b_2} B_5(j) \frac{\Gamma(x_2+1)}{\Gamma{(x_2+B_5(b_2)+1)}} \sum_{k_{b_2} = 0}^{x_2} \frac{\Gamma(k_{b_2}+B_5(b_2))}{\Gamma{(k_{b_2}+B_5(b_2-1)+1)}}\cdots \sum_{k_1 = 0}^{k_2} B_5(1) \frac{\Gamma(k_{1}+B_5(1))}{\Gamma(k_{1}+1)} k_1^{\delta_4}\\
      & \sssim{e_2} \frac{\prod_{j = 1}^{b_2} B_5(j)}{x_2^{B_5(b_2)}} \int_0^{x_2} \dd k_{b_2} k_{b_2}^{B_5(b_2)-B_5(b_2-1)-1} \int_0^{k_{b_2}} \cdots k_2^{B_5(2)-B_5(1)-1}\int_0^{k_2} \dd k_1\, k_1^{B_5(1)-1+\delta_4 }\\& \sssim{e_2} \frac{\prod_{j = 1}^{b_2}B_5(j)}{\prod_{j = 1}^{b_2}(B_5(j)+\delta_4)} x_2^{\delta_4} = \frac{x_2^{\delta_4}}{\prod_{j = 1}^{b_2}(1+\delta_4/B_5(j))} ~.\label{e:productbound}
     \end{equs}
     We note that the product on the right hand side of the above expression converges for large values of $b_2$.

      We now proceed to evaluate the integral term of \eref{e:prop}. Choosing $h_4$ as in \eref{e:tt4} and summing over all paths $\gamma$ from $(b_2,x_2) \to \TT_4^*$ we obtain
      we have
     \begin{equs}
      \Ex{(b_2,x_2)}{\int_0^\tau h_4(X_t) \dt } &= h_4 \sum_{\gamma~:~(b_2,x_2) \to \TT_4^*} \px{(b_2,x_2)}{\gamma} \sum_{x \in \gamma} x_1^{\delta_4'}x_2^{\delta_4''} \Ex{x}{\tau_x}\\&= h_4 \sum_{k_{b_2}= 0}^{x_2} \px{{(b_2,x_2)}}{\leftarrow_{(b_2,k_{b_2})}} \Bigg[b_2^{\delta_4'}\pc{ \frac{k_{b_2}^{\delta_4''}}{k_{b_2}^3 + k_{b_2}^2 B_5(b_2)} + \sum_{j = k_{b_2}}^{x_2} \frac{j^{\delta_4''}}{j^3 + j^2 B_5(b_2)}} +  \\& \qquad  \cdots + \sum_{k_{i}= 0}^{k_{i+1}} \px{(i,k_{i+1})}{\leftarrow_{(i,k_{i})}} \Bigg[i^{\delta_4'}\pc{\frac{k_{i}^{\delta_4''}}{k_{i}^3 + k_{i}^2 B_5(i)} + \sum_{j = k_{i}}^{k_{i+1}} \frac{j^{\delta_4''}}{j^3 + j^2 B_5(i)}} + \\&\qquad \qquad \cdots + \sum_{k_{1}= 0}^{k_2} \px{{(1,k_2)}}{\leftarrow_{(1,k_{1})}} \pc{\frac{k_{1}^{\delta_4''}}{k_{1}^3 + k_{1}^2 B_5(1)} + \sum_{j = k_{1}}^{k_2} \frac{j^{\delta_4''}}{j^3 + j^2 B_5(1)}}\Bigg]\dots \Bigg]\,.
     \end{equs}
     Furthermore, because we are in the limit $x_2 \to \infty$ we approximate the terms in round brackets as
     \[\frac{k_{i}^{\delta_4''}}{k_{i}^3 + k_{i}^2 B_5(i)} + \sum_{j = k_{i}}^{k_{i+1}} \frac{j^{\delta_4''}}{j^3 + j^2 B_5(i)} \sssim{e_2} \int_{k_i}^{k_{i+1}}\dd j\, j^{\delta_4''-3} = \frac{1}{\delta_4''-2} \pc{k_{i+1}^{\delta_4''-2}- k_i^{\delta_4''-2}}\,,\]
     we obtain
     \begin{equs}
      \Ex{(b_2,x_2)}{\int_0^\tau h(X_t) \dt } & \sssim{e_2} h_4 \sum_{i = 1}^{b_2} i^{\delta_4'} \frac{\prod_{j = i}^{b_2}B_5(j)}{x_2^{B_5(b_2)}} \int_0^{x_2} \dd k_{b_2} k_{b_2}^{B_5(b_2)-B_5(b_2-1)-1} \\& \qquad \qquad \qquad \qquad\cdots k_{i+1}^{B_5(i+1)-B_5(i)-1}\int_0^{k_i+1} \dd k_i\,k_{i}^{B_5(i)-1} \pc{k_{i+1}^{\delta_4''-2} - k_{i}^{\delta_4''-2}}\\
      &\sssim{e_2} h_4 x_2^{\delta_4'' - 2}\sum_{i = 1}^{b_2} \frac{i^{\delta_4'}}{(\delta_4''-2) + B_5(i)} \prod_{j = i+1}^{b_2} \frac{1}{1 + (\delta_4''-2)/B_5(j)} ~.
     \end{equs}
     Recalling that the product \eref{e:productbound} converges for large values of $b_2$, we approximate its value by a constant that can be absorbed by $m_4^*$. Through a similar approximation, we write the right hand side of the expression above, for large values of $b_2$, as $h_4 x_2^{\delta_4'' - 2}\sum_{i = 1}^{b_2} ((\delta_4''-2) + B_5(i))^{-1}$. Combining these two approximations we obtain our candidate for the local Lyapunov function in $\TT_4$:
     \begin{equ}
      V_4(x) = {m_4} x_2^{\delta_4} +  h_4 x_2^{\delta_4'' - 2}\sum_{i = 1}^{b_2} \frac{i^{\delta_4'}}{(\delta_4''-2) + B_5(i)} ~.
     \end{equ}
     This function exists for all $\delta_4' \in \Rr$, $\delta_4, \delta_4'' \in \RR$ and scales homogeneously under $\slw$ for $w=e_2$ if $\delta_4'' = \delta_4 + 2$.
    \end{proof}

\section{Alternative scaling procedure}

In this section we present an alternative method for the construction of the Lyapunov funtion in the transport regions $\TT_1, \TT_2, \TT_3$. In Section~\ref{s:transport} we required that the candidate Lyapunov function $V_i$ scales homogeneously under \emph{all} scaling transformations mapping certain subsets of $\WW_i$ onto themselves. Now we simply require for $V_i$ to scale homogeneously under the scaling transformation that maps $\TT_i$ onto a set of $\RR^2$ that includes $\TT_i$ itself. In other words, we need the preimage of $\TT_i$ under $\slw$ to be in $\TT_i$ itself.
This way, by the homogeneous scaling assumption it is sufficient to define $V_i$ on a compact subset of $\TT_i$ to know $V_i$ in all $\TT_i$. Note that for all $\TT_1, \TT_2, \TT_3$ we have that such a scaling is given by $\slw$ with $w = (1,1)/\sqrt 2$.

\lmm{ \label{l:v3a} For an initial condition $x_0 =  (x_1, x_2) \in \TT_3$, the Lyapunov function $\obar V_3$ solving \eref{e:bvp} with
\begin{equ} \label{e:tt3a}
h_3(x) := h_3 x_2^{\delta_3''} \lambda_3(x) \qquad \text{and} \qquad V_4(x) = m_4 x_2^{\delta_4}\,,
\end{equ}
for $h_3 > 0$ and $m_4 = m_4(b_2) > 0$ is well defined for all $\delta_3'' \in \Rr$ and $ \delta_4> 0$. Furthermore, for the choice of constants
\begin{equ}\label{e:hom3a}
  \delta_3'' = \delta_4\,,
\end{equ}
we can write $\obar V_3 (x) = x_2^{\delta_4} (m_4 + h_3(x_1-b_2))$.
}
\begin{proof}
  By the method of characteristics we obtain $\obar V_3$ satisfying \eref{e:bvp} by integrating $h_3$ along the solutions of the set of ordinary differential equations $\dot x = T_3 x$. Recalling by \eref{e:transportgen} that such solutions are moving from $x_1(0)$ to $b_2$ on lines with $x_2(t) = x_2(0)$, by our choice \eref{e:tt3a} of boundary condition on $\TT_{34}$ we obtain
  \begin{equ}
   \obar V_3(x) = m_4 x_2^{\delta_4} + \int_{b_2}^{x_1} h_3 x_2^{\delta_3''} \lambda_3((z,x_2))  \frac{1}{\lambda_3((z,x_2))} \dd z =  m_4 x_2^{\delta_4} - { h_3 } b_2 x_2^{\delta_3''} + { h_3 } x_1 x_2^{\delta_3''}   ~.
  \end{equ}
  This function is clearly well defined in $\RR^2$ for all choices of parameters. Now we see that in order for $\obar V_3$ to scale homogeneously under $\slw$ for all $w \in \WW_3$ we need to have \eref{e:hom3a} for $\delta_3' > 4$ as $h_3 > 0$. This directly implies that $\obar V_3$ has the desired form.
\end{proof}

We note that the function $\obar V_3$ defined in \lref{l:v3a} does not scale homogeneously under $\slw$. For this reason, proceed as in \cite{mattingly151} and introduce a new ``dummy'' coordinate $\lambda$ and a scaling transformation
\begin{equ}
  \mathscr S_l^{(1,1,1)}~:~(x_1,x_2,\chi) \mapsto (l x_1, l x_2, l \chi)\,.
\end{equ}
Then, we define the Lyapunov function $\obar V_3$ in the new set of coordinates as
\begin{equ}
  \obar V_3((x_1,x_2,\chi)) := x_2^{\delta_4} (\chi m_4 + h_3(x_1- \chi b_2))\,.
\end{equ}
It is manifest that this function scales homogeneously under $\mathscr S_l^{(1,1,1)}$.

\lmm{ \label{l:v2a} For an initial condition $x_0 = (x_1, x_2) \in \TT_2$, the Lyapunov function $\obar V_2$ solving \eref{e:bvp} with
\begin{equ}\label{e:tt2a}
h_2(x) := h_2 ( x_1 + 5 x_2)^{\delta_2'} \lambda_3(x) \qquad \text{and} \qquad V_3(\pi_{23} x) = m_3 (x_1 + 5 x_2)^{\delta_3}\,,
\end{equ}
for $h_2 > 0$ is well defined for all $\delta_2' \in \Rr$ and $ \delta_3> 0$. Furthermore, for the choice of constants
\begin{equ}
  \delta_2' = \delta_3''\,,
\end{equ}
we have $\obar V_2 \ssim{w} l^{\delta_2}$ for $w = (1,1,1)$ and $\delta_2 := \delta_2'+1$.}
\begin{proof}
We again find the solution to \eref{e:bvp} in $\TT_2$ by the method of characteristics. Denoting by $\gamma_3(x,\pi_{23})$ the path along the characteristic of $T_2$ starting at $x$ and ending at $\pi_{23}x$ and noting that $h_2(\,\cdot\,)$ defined in \eref{e:tt2a} is constant on such a path we have
\begin{equ}
 \obar V_2(x) = \obar V_3(\pi_{23}x) + h_2(x_1+5x_2)^{\delta_2'} \int_{\gamma_3(x,\pi_{23}x)} \lambda_3(z) \frac{1}{\lambda_3(z)} \dd z\,.
\end{equ}
Consequently, using that $\pi_{23} x = (x_1+5x_2)(1+5 b_1)^{-1}(1,b_1)$ we write the explicit result of the integral as
\begin{equ}\label{e:explicitv2a}
  \obar V_2(x) = (x_1+5x_2)^{\delta_4} \pc{ \chi m_4 + h_3\pc{\frac{x_1+5x_2}{1+5 b_1}- \chi b_2}} + {h_2}(x_1+5x_2)^{\delta_2'}  \mathscr L_{b_1}(x)\,,
\end{equ}
for $\mathscr L_{b_1}(x) := x_1 \sqrt{1 + b_1^{-2}} \sin(\arctan(b_1^{-1}) - \arctan(1/5) )/\sin(\arctan(x_2/x_1)+\arctan(b_1^{-1}))$. By the homogeneous scaling property of $\obar V_3$ under $\mathscr S_l^{(1,1,1)}$ we have that $\obar V_3(\pi_{23} x) = m_3(x_1 + 5 x_2)^{\delta_3}$ for  $\delta_3 := \delta_3''+1$ and $m_3 = m_3(b_1) := h_3(b_1)^{\delta_3''}(1+b_1)^{\delta_3}$.
It is easy to see that the function $\obar V_2$  from \eref{e:explicitv2a} scales homogeneously under $\mathscr S_l^{(1,1,1)}$.
\end{proof}

\lmm{\label{l:v1a}For an initial condition $x_0 = (x_1, x_2) \in \TT_1$, the function $\obar V_1$ solving \eref{e:bvp} with
\begin{equ}
h_1(x) := h_1 x_1^{\delta_1'}\lambda_{3}(x) \qquad \text{and} \qquad \obar V_2(x) = m_2 x_1^{\delta_2}\,,
\end{equ}
for $h_1 > 0$ is well defined for all $\delta_1' \in \Rr$ and $ \delta_2> 0$. Furthermore, for the choice of constants
\begin{equ}\label{e:hom1a}
  \delta_2' = \delta_1'\,,
\end{equ}
we can write $\obar V_1 (x) = x_1^{\delta_2'} (m_2 x_1 - {h_1 } { (x_2-x_1/b_1) })$.
}
\begin{proof}
  We obtain the Lyapunov function by integrating along the characteristic lines of the transport operator $T_1$. Noting that these lines satisfy $x_1(t) = x_1(0)$ for all $t > 0$ we write
  \begin{equ}
    \obar V_1(x) = \obar V_2(\pi_{12} (x)) + h_1 \int_{x_2}^{x_1/b_1} x_1^{\delta_1'} \lambda_3((x_1,y)) \frac{1}{\lambda_3((x_1,y))} \dd y = m_2(x) x_1^{\delta_2} - {h_1 } x_1^{\delta_1'} { (x_2-x_1/b_1) } \,, \label{e:v1a}
  \end{equ}
  for $\pi_{12} (x) = (x_1,x_1/b_1)$ and
  \[m_2(x) := x_1^{\delta_2}(1+5/b_1)^{\delta_2} \pc{\chi m_4 + h_3\pc{x_1+\frac{1+5/b_1}{1+5 b_1} - \chi b_2}} + h_2 x_1^{\delta_2'} (1+b_1/5)^{\delta_2'}\,.\]
   We immediately recognize that the right hand side of \eref{e:v1a} scales homogeneously under $\mathscr S_l^{(1,1,1)}$ if \eref{e:hom1a} holds, resulting in the desired definition for $\obar V_1(\,\cdot\,)$\,.
\end{proof}

\end{document}